\documentclass[12pt]{article}%
\usepackage{amsmath}
\usepackage{amsfonts}
\usepackage{amssymb}
\usepackage{graphicx}
\usepackage{color}%
\providecommand{\U}[1]{\protect\rule{.1in}{.1in}}
\newtheorem{theorem}{Theorem}[section]

\newtheorem{conjecture}{Conjecture}
\newtheorem{corollary}[theorem]{Corollary}

\newtheorem{lemma}[theorem]{Lemma}

\newtheorem{observation}[theorem]{Observation}
\newtheorem{problem}{Problem}
\newtheorem{proposition}[theorem]{Proposition}

\providecommand{\boksie}{\ensuremath{\mathbin{\raisebox{0.3mm}{$\scriptstyle\square$}}}}
\begin{document}

\title{\textbf{Irredundance Graphs}}
\author{C.M. Mynhardt\thanks{Supported by the Natural Sciences and Engineering
Research Council of Canada.}\\Department of Mathematics and Statistics\\University of Victoria, Victoria, BC, \textsc{Canada}\\{\small kieka@uvic.ca}
\and A. Roux\\Department of Mathematical Sciences\\Stellenbosch University, Stellenbosch, \textsc{South Africa}\\{\small rianaroux@sun.ac.za}}
\maketitle

\begin{abstract}
A set $D$ of vertices of a graph $G=(V,E)$ is irredundant if each $v\in D$
satisfies (a) $v$ is isolated in the subgraph induced by $D$, or (b) $v$ is
adjacent to a vertex in $V-D$ that is nonadjacent to all other vertices in
$D$. The upper irredundance number $\operatorname{IR}(G)$ is the largest
cardinality of an irredundant set of $G$; an $\operatorname{IR}(G)$-set is an
irredundant set of cardinality $\operatorname{IR}(G)$.

The $\operatorname{IR}$-graph of $G$ has the irredundant sets of $G$ of
maximum cardinality, that is, the $\operatorname{IR}(G)$-sets, as vertex set,
and sets $D$ and $D^{\prime}$ are adjacent if and only if $D^{\prime}$ is
obtained from $D$ by exchanging a single vertex of $D$ for an adjacent vertex
in $D^{\prime}$. We study the realizability of graphs as $\operatorname{IR}%
$-graphs and show that all disconnected graphs are $\operatorname{IR}$-graphs,
but some connected graphs (e.g. stars $K_{1,n},\ n\geq2$, $P_{4}%
,\ P_{5},\ C_{5},\ C_{6},\ C_{7}$) are not. We show that the double star
$S(2,2)$ -- the tree obtained by joining the two central vertices of two
disjoint copies of $P_{3}$ -- is the unique smallest $\operatorname{IR}$-tree
with diameter $3$ and also a smallest non-complete $\operatorname{IR}$-tree,
and the tree obtained by subdividing a single pendant edge of $S(2,2)$ is the
unique smallest $\operatorname{IR}$-tree with diameter $4$.

\end{abstract}

\noindent\textbf{Keywords:\hspace{0.1in}}Irredundance; Reconfiguration
problem; $\operatorname{IR}$-graph

\noindent\textbf{AMS Subject Classification Number 2010:\hspace{0.1in}}05C69

\section{Introduction}

\label{Sec_Intro}Reconfiguration problems are concerned with determining
conditions under which a feasible solution to a given problem can be
transformed into another such solution via a sequence of feasible solutions in
such a way that any two consecutive solutions are adjacent according to a
specified adjacency relation. The solutions and their adjacencies form the
vertex and edge sets, respectively, of the associated reconfiguration graph.
Typical questions about reconfiguration graphs concern their structure
(connectedness, Hamiltonicity, diameter, planarity), realizability (which
graphs can be realized as a specific type of reconfiguration graph), and
algorithmic properties (finding shortest paths between solutions quickly). 

Reconfiguration versions of graph colouring and other graph problems, such as
independent sets, cliques, and vertex covers, have been studied in
e.g.~\cite{Beier, BC2, CHJ1, CHJ2, Ito2, IKD}.

Domination reconfiguration problems involving (not necessarily minimal)
dominating sets of different cardinalities were first considered in 2014 by
Haas and Seyffarth \cite{HS1} and subsequently also in, for example,
\cite{davood, HS2, haddadan-2015, MRT}. Domination reconfiguration problems
involving only minimum-cardinality dominating sets were introduced by Fricke,
Hedetniemi, Hedetniemi, and Hutson \cite{FHHH11}, and also studied in
\cite{Bien, CHH10, Dyck, EdThesis, Laksh, Laura1, Amu, SS08, SMN}. We study
the upper irredundance graph ($\operatorname{IR}$-graph for short) of a given
graph $G$ -- the ways in which maximum irredundant sets (defined below) of $G$
can be reconfigured successively into other such sets by exchanging (swapping)
a single vertex for a neighbour in each step.

As one of our main results we show in Section \ref{Sec_Disc} that all
disconnected graphs are $\operatorname{IR}$-graphs (Theorem \ref{Thm_Disconn}%
). In Proposition \ref{Prop_diam2} we show that every $\operatorname{IR}%
$-graph with diameter $2$ contains an induced $4$-cycle; this result implies
that not all graphs can be realized as $\operatorname{IR}$-graphs. For
example, the stars $K_{1,n}$ are not $\operatorname{IR}$-graphs. Generalizing
this result we show in Proposition \ref{Prop_K1_n} that non-complete graphs
with universal vertices are also not $\operatorname{IR}$-graphs. Our other
main results are given in Section \ref{Sec_Trees}, where we characterize the
smallest $\operatorname{IR}$-trees, that is, trees that are $\operatorname{IR}%
$-graphs, of diameter $3$ or $4$. The aforementioned results culminate in
Theorem \ref{Thm_Main}, which states that the cycles $C_{5},C_{6},C_{7}$ and
the paths $P_{3},P_{4},P_{5}$ are not $\operatorname{IR}$-graphs, the only
connected $\operatorname{IR}$-graphs of order $4$ are $K_{4}$ and $C_{4}$, and
the double star $S(2,2)$ is a the smallest non-complete $\operatorname{IR}%
$-tree. In the rest of Section \ref{Sec_Intro} we give basic definitions
(Section \ref{Sec_Defs}) and the definition of an $\operatorname{IR}$-graph,
first defined by Fricke et al. \cite{FHHH11} (Section \ref{Sec_IRGr}). Section
\ref{Sec_Basic} contains elementary results and lemmas required later.
Problems for future research are given in Section \ref{Sec_Qs}.

\subsection{Definitions}

\label{Sec_Defs}In general, we follow the notation of \cite{CLZ}. For
domination related concepts not defined here we refer the reader to
\cite{HHS}. Given a graph $G=(V,E)$, the \emph{open }and \emph{closed
neighbourhoods of a vertex }$v$ of $G$ are, respectively, $N(v)=\{u\in V:uv\in
E\}$ and $N[v]=N(v)\cup\{v\}$. The \emph{open }and \emph{closed }%
neighbourhoods of a set $D\subseteq V$ are, respectively, $N(D)=\bigcup_{v\in
D}N(v)$ and $N[D]=N(D)\cup D$. The set $D$ \emph{dominates} a set $A\subseteq
V(G)$ if $A\subseteq N[D]$, and is a \emph{dominating set }of $G$ if
$N[D]=V(G)$.

A $D$-\emph{private neighbour }of $v\in D$ is a vertex $v^{\prime}$ that is
dominated by $v$ but by no vertex in $D-\{v\}$. The set of $D$-private
neighbours of $v$ is called the \emph{private neighbourhood of }$v$ \emph{with
respect to }$D$ and denoted by $\operatorname{PN}(v,D)$, that is,
$\operatorname{PN}(v,D)=N[v]-N[D-\{v\}]$. A dominating set $D$ is a
\emph{minimal dominating set }if no proper subset of $D$ is a dominating set.
It is well known \cite[Theorem 1.1]{HHS} that a dominating set $D$ is minimal
dominating if and only if each $v\in D$ has a $D$-private neighbour. The
\emph{lower and upper domination numbers} of $G$ are the cardinalities of a
smallest dominating set and a largest minimal dominating set, respectively,
and are denoted by $\gamma(G)$ and $\Gamma(G)$, respectively. A $\gamma
$-\emph{set }is a dominating set of cardinality $\gamma(G)$, and a $\Gamma
$-\emph{set }is a minimal dominating set of cardinality$~\Gamma(G)$.

The concept of irredundance was introduced by Cockayne, Hedetniemi and Miller
\cite{CHM} in 1978. A set $D\subseteq V$ is \emph{irredundant }if
$\operatorname{PN}(v,D)\neq\varnothing$ for each $v\in D$, and \emph{maximal
irredundant} if no superset of $D$ is irredundant. The \emph{irredundance
number }$\operatorname{ir}(G)$ is the minimum cardinality of a maximal
irredundant set of $G$, and the \emph{upper irredundance number }%
$\operatorname{IR}(G)$ is the largest cardinality of an irredundant set of
$G$. An $\operatorname{ir}$-\emph{set} of $G$ is a maximal irredundant set of
cardinality $\operatorname{ir}(G)$, and an $\operatorname{IR}$-\emph{set} of
$G$, sometimes also called an $\operatorname{IR}(G)$-\emph{set}, is an
irredundant set of cardinality$~\operatorname{IR}(G)$.

\begin{figure}[ptb]%
\centering
\includegraphics[
height=1.6777in,
width=5.5011in
]%
{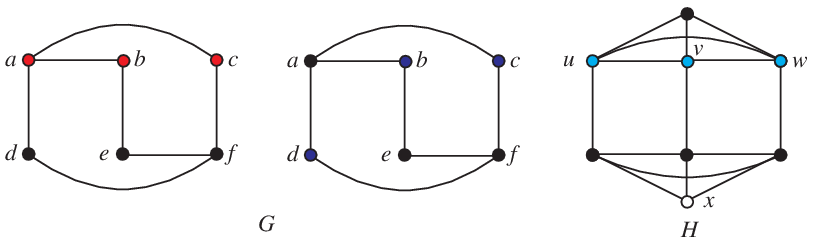}%
\caption{A graph $G$ and $\operatorname{IR}(G)$-sets
$A=\{a,b,c\},\ B=\{b,c,d\}$, and a graph $H$ with non-dominating
$\operatorname{IR}$-set $\{u,v,w\}$}%
\label{Fig_IRset}%
\end{figure}

Let $D$ be an irredundant set of $G$. For $v\in D$, it is possible that
$v\in\operatorname{PN}(v,D)$; this happens if and only if $v$ is isolated in
the subgraph $G[D]$ induced by $D$. If $u\in\operatorname{PN}(v,D)$ and $u\neq
v$, then $u\in V-D$; in this case $u$ is an \emph{external }$D$-\emph{private
neighbour of }$v$. The set of external $D$-private neighbours of $v$ is
denoted by $\operatorname{EPN}(v,D)$. An isolated vertex of $G[D]$ may or may
not have external $D$-private neighbours, but if $v$ has positive degree in
$G[D]$, then $\operatorname{EPN}(v,D)\neq\varnothing$. Figure \ref{Fig_IRset}
shows a graph $G$ and two of its $\operatorname{IR}$-sets $A=\{a,b,c\}$ and
$B=\{b,c,d\}$. Each vertex in $A$ has positive degree in $G[A]$, and
$\operatorname{PN}(a,A)=\operatorname{EPN}(a,A)=\{d\}$, $\operatorname{PN}%
(b,A)=\operatorname{EPN}(b,A)=\{e\}$ and $\operatorname{PN}%
(c,A)=\operatorname{EPN}(c,A)=\{f\}$. On the other hand, $B$ is an independent
set; $\operatorname{PN}(b,B)=\{b,e\}$ and $\operatorname{EPN}(b,B)=\{e\}$,
$\operatorname{PN}(c,D)=\{c\}$, $\operatorname{PN}(d,D)=\{d\}$ and
$\operatorname{EPN}(c,D)=\varnothing=\operatorname{EPN}(d,D)$. Dominating,
independent and irredundant sets are related as follows (see
\cite[Propositions 3.8 and 3.9]{HHS}).

\begin{observation}
\label{Ob_irr_dom}

\begin{enumerate}
\item[$(i)$] If a set is irredundant and dominating, it is maximal irredundant
and minimal dominating.

\item[$(ii)$] A dominating set is minimal dominating if and only if it is irredundant.

\item[$(iii)$] Any independent set is irredundant, and any maximal independent
set is minimal dominating and maximal irredundant.
\end{enumerate}
\end{observation}

However, a maximal irredundant set need not be dominating; the
$\operatorname{IR}$-set $\{u,v,w\}$ of the graph $H$ in Figure \ref{Fig_IRset}
is an example of such a set -- the addition of any vertex to dominate $x$ will
destroy the private neighbours of the other vertices.
%TCIMACRO{\FRAME{ftbpFU}{5.5011in}{1.6777in}{0pt}{\Qcb{A graph $G$ and
%$\operatorname{IR}(G)$-sets $A=\{a,b,c\},\ B=\{b,c,d\}$, and a graph $H$ with
%non-dominating $\operatorname{IR}$-set $\{u,v,w\}$}}{\Qlb{Fig_IRset}%
%}{ir_set.eps}{\special{ language "Scientific Word";  type "GRAPHIC";
%maintain-aspect-ratio TRUE;  display "USEDEF";  valid_file "F";
%width 5.5011in;  height 1.6777in;  depth 0pt;  original-width 5.444in;
%original-height 1.6414in;  cropleft "0";  croptop "1";  cropright "1";
%cropbottom "0";  filename 'IR_set.eps';file-properties "XNPEU";}} }%
%BeginExpansion
%EndExpansion

The complete bipartite graph $K_{1,k},\ k\geq1$, is called a \emph{star}. Let
$K_{1,k}$ have partite sets $\{u\}$ and $\{v_{1},...,v_{k}\}$. The
\emph{(generalized) spider }$\mathrm{Sp}(\ell_{1},...,\ell_{k}),\ \ell_{i}%
\geq1,\,k\geq2$, is a tree obtained from $K_{1,k}$ by subdividing the edge
$uv_{i}$ $\ell_{i}-1$ times, $i=1,...,k$. The \emph{double star }$S(k,n)$ is
the tree obtained by joining the centres of the stars $K_{1,k}$ and $K_{1,n}$.
The \emph{double spider }$\mathrm{Sp}(\ell_{1},...,\ell_{k};m_{1},...,m_{n})$
is obtained from $S(k,n)$ by subdividing the edges of the $K_{1,k}$-subgraph
$\ell_{i}-1$ times, $i=1,...,k$, and the edges of the $K_{1,n}$-subgraph
$m_{i}-1$ times, $i=1,...,n$.

\subsection{$\operatorname{IR}$-Graphs}

\label{Sec_IRGr}First defined by Fricke et al. \cite{FHHH11}, the
\emph{$\gamma$-graph }$G(\gamma)$ of \emph{$G$} has the $\gamma$-sets of $G$
as its vertex set, where $S$ and $S^{\prime}$ are adjacent in $G(\gamma)$ if
and only if there exist vertices $u\in S$ and $v\in S^{\prime}$ such that
$uv\in E(G)$ and $S^{\prime}=(S-\{u\})\cup\{v\}$. This model of adjacency in
$G(\gamma)$ is referred to as the \emph{slide-model}; other authors, for
example \cite{SS08} and especially \cite{Dyck}, consider a \emph{jump-model},
where the vertices $u^{\prime}$ and $v^{\prime}$ need not be adjacent. Fricke
et al. showed that every tree is the $\gamma$-graph of some graph and
conjectured that every graph is the $\gamma$-graph of some graph; Connelly et
al. \cite{CHH10} proved this conjecture to be true. On the other hand, not all
graphs are $\gamma$-graphs if one uses the jump-model; \textquotedblleft
jump\textquotedblright\ $\gamma$-graphs were characterized in \cite{Dyck}. For
additional results on $\gamma$-graphs, see \cite{CHH10, EdThesis, FHHH11}.

As observed in \cite{Laura1}, the construction in \cite{CHH10} also suffices
to prove that every graph is the $\pi$-graph (according to the slide-model) of
infinitely many graphs, where $\pi$ is any of the parameters $\gamma$,
$\operatorname{ir}$, $\gamma_{\operatorname{pr}}$ (the paired-domination
number), $\gamma_{t}$ (the total domination number), and $\gamma_{c}$ (the
connected domination number). Different constructions in \cite{Laura1} further
show that every graph is the $\pi$-graph (again according to the slide-model)
of infinitely many graphs, for a variety of domination related parameters
$\pi$, including $\Gamma$. The study of $\operatorname{IR}$-graphs is
mentioned as an open problem in \cite{Laura1}; we initiate the study of these
graphs here.

Following \cite{FHHH11}, we define the \emph{$\operatorname{IR}$-graph
}$G(\operatorname{IR})$\emph{ }of \emph{$G$} to be the graph whose vertex set
consists of the $\operatorname{IR}(G)$-sets, i.e., the irredundant sets of $G$
of maximum cardinality, where two $\operatorname{IR}(G)$-sets $D$ and
$D^{\prime}$ are adjacent if and only if there exist vertices $u\in D$ and
$v\in D^{\prime}$ such that $uv\in E(G)$ and $D^{\prime}=(D-\{u\})\cup\{v\}$.
We shorten the expression $D^{\prime}=(D-\{u\})\cup\{v\}$ to
$D\overset{uv}{\sim}D^{\prime}$, and also write $D\sim_{H}D^{\prime}$ to show
that $D$ and $D^{\prime}$ are adjacent in $H=G(\operatorname{IR})$. When
$D\overset{uv}{\sim}D^{\prime}$, we say that $v$ is \emph{swapped into} and
$u$ is \emph{swapped out of} the $\operatorname{IR}(G)$-set, or simply that
$u$ and $v$ are \emph{swapped}. The notation $u\sim v$ ($u\nsim v$,
respectively) indicates that $u$ is adjacent (respectively nonadjacent) to
$v$; we sometimes write $u\sim_{G}v$ for emphasis.

To prove that a given graph $H$ is an $\operatorname{IR}$-graph, one needs to
construct a graph $G$ such that $G(\operatorname{IR})\cong H$. In this
situation the graphs $G$ and $H$ are referred to as the \emph{source} and
\emph{target }graphs, respectively.

\section{Basic Results}

\label{Sec_Basic}We begin by stating a few simple results on
$\operatorname{IR}$-graphs. Following \cite{CLZ} we denote the (disjoint)
union of graphs $G_{1}$ and $G_{2}$ by $G_{1}+G_{2}$, and their Cartesian
product by $G_{1}\boksie G_{2}$. If $G_{i}\cong G$ for $i=1,...,n$, the union
$G_{1}+\cdots+G_{n}$ is denoted by $nG$ and the Cartesian product
$G_{1}\boksie\cdots\boksie G_{n}$ by $G^{n}$.

Clearly, $\operatorname{IR}(G)=1$ if and only if $G=K_{n},\ n\geq1$. Hence
$K_{n}$ has $n$ $\operatorname{IR}$-sets, and any two of them are adjacent in
$K_{n}(\operatorname{IR})$, that is, $K_{n}(\operatorname{IR})=K_{n}$.

\begin{proposition}
\label{Prop_Qn}

\begin{enumerate}
\item[$(i)$] If $H_{1}$ and $H_{2}$ are $\operatorname{IR}$-graphs, then
$H_{1}\boksie H_{2}$ is an $\operatorname{IR}$-graph.

\item[$(ii)$] For all $n\geq1$, the hypercube $Q_{n}$ is an $\operatorname{IR}%
$-graph. In particular, $C_{4}\cong Q_{2}$ is an $\operatorname{IR}$-graph.
\end{enumerate}
\end{proposition}

\noindent\textbf{Proof.\hspace{0.1in}}$(i)$ If $H_{1}=G_{1}(\operatorname{IR}%
)$ and $H_{2}=G_{2}(\operatorname{IR})$, then $H_{1}\boksie H_{2}=(G_{1}%
+G_{2})(\operatorname{IR})$.

\noindent$(ii)$ Let $G=nK_{2}$. Then $G(\operatorname{IR})=(K_{2})^{n}\cong
Q_{n}$.~$\blacksquare$

\bigskip

The next result is of independent interest, and is also used throughout the
paper to find more $\operatorname{IR}$-sets by using external private
neighbours in a given $\operatorname{IR}$-set.\ For an irredundant set $X$, we
weakly partition $X$ into subsets $Y$ and $Z$ (one of which may be empty),
where each vertex in $Z$ is isolated in $G[X]$ and each vertex in $Y$ has at
least one external private neighbour. (This partition is not necessarily
unique. Isolated vertices of $G[X]$ with external private neighbours can be
allocated arbitrarily to $Y$ or $Z$.) For each $y\in Y$, let $y^{\prime}%
\in\operatorname{EPN}(y,X)$ and define $Y^{\prime}=\{y^{\prime}:y\in Y\}$. Let
$X^{\prime}=(X-Y)\cup\{Y^{\prime}\}$; note that $|X|=|X^{\prime}|$. We call
$X^{\prime}$ the \emph{flip-set }of $X$, or to be more precise, the
\emph{flip-set }of $X$ \emph{using }$Y^{\prime}$.

\begin{proposition}
\label{Prop_IR_sets}If $X$ is an $\operatorname{IR}$-set of $G$, then so is
any flip-set $X^{\prime}$ of $X$.
\end{proposition}

\noindent\textbf{Proof.\hspace{0.1in}}Consider any $x\in X^{\prime}$. We show
that $x$ has an $X^{\prime}$-private neighbour. With notation as above,
suppose $x\in Z=X-Y=X^{\prime}-Y^{\prime}$. By definition of $Z$, $x$ is
isolated in $G[X]$. Since each vertex in $Y^{\prime}$ is an $X$-external
private neighbour of some $y\in Y$, no vertex in $Y^{\prime}$ is adjacent to
$x$. Therefore $x$ is isolated in $G[X^{\prime}]$, that is, $x\in
\operatorname{PN}(x,X^{\prime})$.

Hence assume $x\in Y^{\prime}$. Then $x=y^{\prime}$ for some $y\in Y$, so
$y^{\prime}$ is an $X$-external private neighbour of $y\in V(G)-X^{\prime}$.
We show that $y$ is an $X^{\prime}$-external private neighbour of $y^{\prime}%
$. Now $y$ is non-adjacent to all vertices in $Z$ because the latter vertices
are isolated in $G[X]$, and $y$ is nonadjacent to all vertices in $Y^{\prime
}-\{y\}$, because each $v^{\prime}\in Y^{\prime}-\{y^{\prime}\}$ is an
$X$-external private neighbour of some $v\in Y-\{y\}$. Therefore
$y\in\operatorname{EPN}(y^{\prime},X^{\prime})$, that is, $y\in
\operatorname{EPN}(x,X^{\prime})$. It follows that $X^{\prime}$ is
irredundant. Since $|X^{\prime}|=|X|$, $X^{\prime}$ is an $\operatorname{IR}%
$-set of $G$.~$\blacksquare$

\section{Disconnected $\operatorname{IR}$-graphs}

\label{Sec_Disc}In this section we resolve the realizability of disconnected
graphs as $\operatorname{IR}$-graphs: all are $\operatorname{IR}$-graphs. The
main idea of the proof is similar to the proof in \cite{CHH10} that all graphs
are $\gamma$-graphs. To show that a given graph $H$ with $V(H)=\{v_{1}%
,...,v_{n}\}$ is a $\gamma$-graph, the authors construct a supergraph $G$ of
$H$ with $\gamma(G)=2$ in which some vertex $u\in V(G)-V(H)$ belongs to all
$\gamma$-sets, these being precisely the sets $\{u,v_{1}\},...,\{u,v_{n}\}$.

Let $H$ be an arbitrary disconnected graph such that $H_{1}$ is one component
of $H$, and $H_{2}$ is the union of all the other components of $H$. Here we
will construct a supergraph $G$ of $H$ such that the $\operatorname{IR}$-sets
of $G$ are of two types $R_{i}$ and $S_{j}$. Each $R_{i}$ is a set of the form
$\{u_{i}\}\cup X$, where $u_{i}\in V(H_{1})$ and $X\subseteq V(G)-V(H)$. Each
$S_{j}$ is of the form $\{v_{j}\}\cup Y$, where $v_{j}\in V(H_{2})$ and
$Y\subseteq V(G)-V(H)-X$. We will show that the edges between the vertices
$u_{i}$ of $H_{1}$ determine the edges between the vertices $R_{i}$ of
$G(\operatorname{IR})$; similarly, the edges between the vertices $v_{j}$ of
$H_{2}$ determine the edges between the vertices $S_{j}$ of
$G(\operatorname{IR})$. This will show that $G(\operatorname{IR})\cong H$.
Proposition \ref{Prop_IR_sets} explains, to some extent, why the construction
only works for disconnected graphs: the $\operatorname{IR}(G)$-sets used to
construct one component of the target graph are not connected, in
$G(\operatorname{IR})$, to their flip-sets, which are used to construct the
other components.

\begin{theorem}
\label{Thm_Disconn}Every disconnected graph is an $\operatorname{IR}$-graph of
infinitely many graphs $G$.
\end{theorem}

\begin{figure}[ptb]%
\centering
\includegraphics[
height=3.0692in,
width=3.4627in
]%
{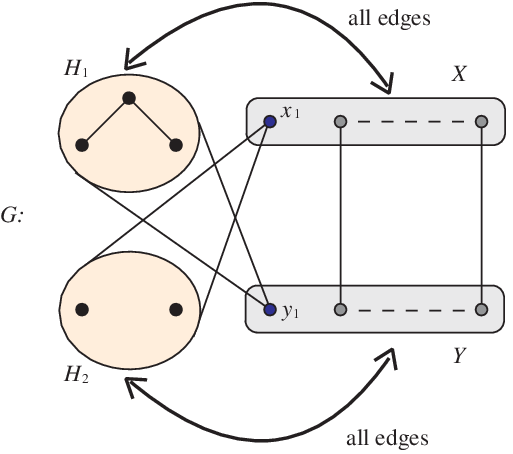}%
\caption{The graph $G$ in the proof of Theorem \ref{Thm_Disconn}}%
\label{Fig_Disconn}%
\end{figure}

\noindent\textbf{Proof.\hspace{0.1in}}Let $H$ be a disconnected graph of order
$n$, let $H_{1}$ be a component of $H$ and let $H_{2}$ be the union of the
other components of $H$. Say $V(H_{1})=\{u_{1},...,u_{n_{1}}\}$ and
$V(H_{2})=\{v_{1},...,v_{n_{2}}\}$. For any $N\geq n$, consider the disjoint
sets $X=\{x_{1},...,x_{N}\}$ and $Y=\{y_{1},...,y_{N}\}$. Construct the
supergraph $G$ of $H$ by adding edges such that $G[X]\cong G[Y]\cong K_{N}$,
every vertex in $X\cup\{y_{1}\}$ is adjacent to every vertex of $H_{1}$, every
vertex in $Y\cup\{x_{1}\}$ is adjacent to every vertex of $H_{2}$, and
$x_{i}\sim y_{i}$ for $i=2,...,N$. See Figure \ref{Fig_Disconn}. Let
$V_{1}=X\cup V(H_{1})$ and $V_{2}=Y\cup V(H_{2})$.%
%TCIMACRO{\FRAME{ftbpFU}{3.4627in}{3.0692in}{0pt}{\Qcb{The graph $G$ in the
%proof of Theorem \ref{Thm_Disconn}}}{\Qlb{Fig_Disconn}}{ir_disconn.eps}%
%{\special{ language "Scientific Word";  type "GRAPHIC";
%maintain-aspect-ratio TRUE;  display "USEDEF";  valid_file "F";
%width 3.4627in;  height 3.0692in;  depth 0pt;  original-width 3.416in;
%original-height 3.0251in;  cropleft "0";  croptop "1";  cropright "1";
%cropbottom "0";  filename 'IR_Disconn.eps';file-properties "XNPEU";}} }%
%BeginExpansion
%EndExpansion

For any $i\in\{1,...,n_{1}\}$, $R_{i}=\{u_{i}\}\cup X$ is an irredundant set
in $G$, because $y_{1}\in\operatorname{EPN}(u_{i},R_{i})$, $v_{1}%
\in\operatorname{EPN}(x_{1},R_{i})$ and $y_{j}\in\operatorname{EPN}%
(x_{j},R_{i})$ for $j\in\{2,...,N\}$. Similarly, for any $i\in\{1,...,n_{2}%
\}$, $S_{i}=\{v_{i}\}\cup Y$ is irredundant. Note that $|R_{i}|=|S_{i}|=N+1$.
We show that $\operatorname{IR}(G)=N+1$ and that the sets $R_{i}$ and $S_{i}$
are the only $\operatorname{IR}$-sets of $G$.

Let $D$ be any maximal irredundant set in $G$. First suppose $D\cap
V(H)=\varnothing$. If, in addition, $D\cap X=\varnothing$ or $D\cap
Y=\varnothing$, then $|D|\leq N$, and if $x_{i},y_{j}\in D$ for some $i,j$,
then $\{x_{i},y_{j}\}$ dominates $G$, hence is maximal irredundant, so that
$D=\{x_{i},y_{j}\}$ and $|D|=2\leq N$.

Now suppose $D$ contains at least two vertices $u_{i}$ of $H_{1}$. (The result
is similar if $D$ contains at least two vertices of $H_{2}$.) If $D$ contains
at least one vertex in $X\cup\{y_{1}\}$, then the $u_{i}$ are not isolated in
$G[D]$ and have no private neighbours in $V_{1}\cup\{y_{1}\}$, and they have
no neighbours in $V_{2}-\{y_{1}\}$. Hence $\operatorname{PN}(u_{i}%
,D)=\varnothing$, a contradiction. Therefore $D\cap(X\cup\{y_{1}%
\})=\varnothing$ and $\operatorname{PN}(u_{i},D)\subseteq V(H_{1})$. It
follows that $I_{1}=D\cap V(H_{1})$ is an irredundant set of $H_{1}$. If
$y_{j}\in D$ for $j\geq2$, then $\operatorname{PN}(y_{j},D)\subseteq
V_{2}-\{y_{1}\}\subseteq N[y_{j}]$. Therefore $D\cap V_{2}=\{y_{j}\}$ and
$|D|\leq\operatorname{IR}(H_{1})+1\leq N$. On the other hand, if $D\cap
Y=\varnothing$ and $I_{2}$ is any irredundant set of $H_{2}$, then $I_{1}\cup
I_{2}$ is maximal irredundant in $G$. In this case $|D|\leq\operatorname{IR}%
(H_{1})+\operatorname{IR}(H_{2})\leq n\leq N$.

Finally, suppose $D$ contains exactly one vertex $u_{i}$ of $H_{1}$. Since
$u_{i}$ dominates $X\cup\{y_{1}\}$, the private neighbours of all other
vertices in $D$ belong to $V(H_{1})\cup V_{2}-\{y_{1}\}$. If $x_{j}\in D$ for
some~$j$, then $\operatorname{PN}(u_{i},D)=\{y_{1}\}$, hence $D\cap
V_{2}=\varnothing$ and $|D|\leq N+1$. Moreover, if $|D|=N+1$, then $D=R_{i}$.
We may therefore assume that $D\cap V_{1}=\{u_{i}\}$. If $y_{1}\in D$, then
$\{u_{i},y_{1}\}$ dominates $G$, hence $D=\{u_{i},y_{1}\}$. If $y_{j}\in D$
for $j\geq2$, then $\{u_{i},y_{j}\}$ dominates $G-V(H_{1})$. Under the
above-mentioned assumption that $D\cap V(H_{1})=\{u_{i}\}$, it follows that
$\{u_{i}\}$ is maximal irredundant in $H_{1}$. Thus $D=\{u_{i},y_{j}\}$. If
$D\cap Y=\varnothing$, then $D=\{u_{i}\}\cup I$, where $I$ is an irredundant
set of $H_{2}$. In all cases, $|D|\leq\operatorname{IR}(H_{1}%
)+\operatorname{IR}(H_{2})\leq n\leq N$.

It follows that $R_{i}=\{u_{i}\}\cup X$, $i\in\{1,...,n_{1}\}$, and
$S_{i}=\{v_{i}\}\cup Y$, $i\in\{1,...,n_{2}\}$, are the only
$\operatorname{IR}$-sets of $G$. Let $H^{\ast}=G(\operatorname{IR})$ with
$V(H^{\ast})=\{R_{i}:i\in\{1,...,n_{1}\}\}\cup\{S_{i}:i\in\{1,...,n_{2}\}\}$.
Note that $R_{i}\sim_{H^{\ast}}R_{j}$ if and only if $u_{i}\sim_{G}u_{j}$ if
and only if $u_{i}\sim_{H_{1}}u_{j}$, and $S_{i}\sim_{H^{\ast}}S_{j}$ if and
only if $v_{i}\sim_{G}v_{j}$ if and only if $v_{i}\sim_{H_{2}}v_{j}$.
Therefore $H^{\ast}[\{R_{1},...,R_{n_{1}}\}]\cong H_{1}$ and $H^{\ast}%
[\{S_{1},...,S_{n_{2}}\}]\cong H_{2}$. Since $|R_{i}-S_{j}|>1$ for all $i$ and
$j$, no $R_{i}$ is adjacent, in $H^{\ast}$, to any $S_{j}$; hence $H\cong
H^{\ast}$.~$\blacksquare$

\section{Connected $\operatorname{IR}$-graphs}

\label{Sec_Con}In the rest of the paper we study connected $\operatorname{IR}%
$-graphs. The results in this section play an important role in Section
\ref{Sec_Trees}, where we investigate trees of small diameter. In the first
lemma we use Proposition \ref{Prop_IR_sets} to further explore the role of
external private neighbours and flip-sets in the construction of
$\operatorname{IR}$-graphs. Combining Lemmas \ref{Lem_Not_independent} and
\ref{Lem_Indep} leads to the result that not all connected graphs are
$\operatorname{IR}$-graphs.

\begin{lemma}
\label{Lem_Not_independent}If $G$ has an $\operatorname{IR}$-set $X$ in which
at least two vertices have external private neighbours and the
$\operatorname{IR}$-graph $H$ of $G$ is connected, then $H$ has an induced
$4$-cycle containing $X$, or $\operatorname{diam}(H)\geq3$ and $d_{H}%
(X,X^{\prime})\geq3$ for any flip-set $X^{\prime}$ of $X$.
\end{lemma}

\noindent\textbf{Proof.\hspace{0.1in}}Let $X=\{x_{1},...,x_{r}\}$ be an
$\operatorname{IR}$-set of $G$ such that, for some $k\geq2$, $x_{1},...,x_{k}$
have external private neighbours, while $x_{k+1},...,x_{r}$ are isolated in
$G[X]$ (and may or may not have external private neighbours). For $i=1,...,k$,
choose $x_{i}^{\prime}\in\operatorname{EPN}(x_{i},X)$ and let $X^{\prime
}=(X-\{x_{1},...,x_{k}\})\cup\{x_{1}^{\prime},...,x_{k}^{\prime}\}$. By
Proposition \ref{Prop_IR_sets}, $X^{\prime}$ is an $\operatorname{IR}(G)$-set.

Assume first that $k=2$ and let $X_{i}=(X-\{x_{i}\})\cup\{x_{i}^{\prime}\}$,
$i=1,2$. Since $x_{i}^{\prime}\in\operatorname{EPN}(x_{i},X)$ and $x_{1}x_{2}$
is the only possible edge of $G[X]$, each $X_{i}$ is independent and therefore
an $\operatorname{IR}(G)$-set. Moreover, $(X\overset{x_{1}x_{1}^{\prime}%
}{\sim}_{H}X_{1}\overset{x_{2}x_{2}^{\prime}}{\sim}_{H}X^{\prime
}\overset{x_{1}^{\prime}x_{1}}{\sim}_{H}X_{2}\overset{x_{2}^{\prime}%
x_{2}}{\sim}_{H}X)$, that is, $H$ has a $4$-cycle $C$. Since $|X-X^{\prime
}|=|X_{1}-X_{2}|=2$, $X\nsim X^{\prime}$ and $X_{1}\nsim X_{2}$, hence $C$ is
an induced $4$-cycle containing $X$.

Now assume $k\geq3$. Then $|X-X^{\prime}|=k\geq3$, and since $H$ is connected,
at least three swaps are required to reconfigure $X$ to $X^{\prime}$.
Therefore $d_{H}(X,X^{\prime})\geq3$ and so $\operatorname{diam}(H)\geq3$, as
asserted.~$\blacksquare$

\bigskip

We present three immediate corollaries to Lemma \ref{Lem_Not_independent}.

\begin{corollary}
\label{Cor_Not_independent}If $G$ has an $\operatorname{IR}$-set that is not
independent and the $\operatorname{IR}$-graph $H$ of $G$ is connected, then
$H$ has an induced $4$-cycle or $\operatorname{diam}(H)\geq3$.
\end{corollary}

Since the $\operatorname{IR}$-sets $X,X_{1},X_{2},X^{\prime}$ in the first
part of the proof of Lemma \ref{Lem_Not_independent} induce a $C_{4}$, we have
the following result.

\begin{corollary}
\label{Cor_C4}If the $\operatorname{IR}$-graph $H$ of $G$ is connected and $G$
has an $\operatorname{IR}$-set $X$ such that

\begin{enumerate}
\item[$(i)$] $G[X]$ has exactly one edge, or

\item[$(ii)$] $X$ is independent but at least two vertices have $X$-external
private neighbours,
\end{enumerate}

then $X$ lies on an induced $4$-cycle of $H$.
\end{corollary}

In the second part of the proof of Lemma \ref{Lem_Not_independent}, if some
$\operatorname{IR}$-set $X$ of $G$ contains $k$ vertices of positive degree in
$G[X]$, then $d_{H}(X,X^{\prime})\geq k$. The next result follows.

\begin{corollary}
\label{Cor_Number_IR-sets1}If the $\operatorname{IR}$-graph $H$ of $G$ is
connected and $G$ has an $\operatorname{IR}$-set $X$ that contains $k\geq3$
vertices of positive degree in $G[X]$, or with $X$-external private
neighbours, then $\operatorname{diam}(H)\geq k$.
\end{corollary}

We next consider graphs whose $\operatorname{IR}$-sets are independent and
whose $\operatorname{IR}$-graphs are connected. Lemma \ref{Lem_Indep}
illustrates that the $\operatorname{IR}$-graphs of these graphs are cyclic.

\begin{lemma}
\label{Lem_Indep}Let $G$ be a graph, all of whose $\operatorname{IR}$-sets are
independent. If the $\operatorname{IR}$-graph $H$ of $G$ is connected and has
order at least three, then $H$ contains a triangle or an induced $C_{4}$.
\end{lemma}

\noindent\textbf{Proof.\hspace{0.1in}}If $\operatorname{IR}(G)=1$, then
$G=K_{n}$ for some $n$, so $H=K_{n}$. Then $n\geq3$ and $H$ has a triangle.
Hence we assume that $\operatorname{IR}(G)\geq2$. Let $X_{0}=\{x_{1}%
,...,x_{r}\}$ be any $\operatorname{IR}(G)$-set and, without loss of
generality, $X_{1}=\{b,x_{2},...,x_{r}\}$, where $b\sim x_{1}$, is an
$\operatorname{IR}(G)$-set such that $X_{0}\sim_{H}X_{1}$. Let $X_{2}$ be any
other $\operatorname{IR}(G)$-set. If $X_{2}=\{c,x_{2},...,x_{r}\}$, then
$b\sim c$, otherwise $(X_{0}-\{x_{1}\})\cup\{b,c\}$ is an independent set
(hence an $\operatorname{IR}$-set) of larger cardinality than $X_{0}$, which
is impossible. Similarly, $x_{1}\sim c$. But then $X_{0}\overset{x_{1}b}{\sim
}_{H}X_{1}\overset{bc}{\sim}_{H}X_{2}\overset{cx_{1}}{\sim}_{H}X_{0}$ and $H$
has a triangle.

Hence, without loss of generality, $X_{2}=\{b,c,x_{3},...,x_{r}\}$, where
$c\sim x_{2}$. Since $X_{1}$ is independent, $b\nsim x_{2},...,x_{r}$, and
since $X_{2}$ is independent, $c\nsim b,x_{3},...,x_{r}$. Let $X_{3}%
=\{x_{1},c,x_{3},...,x_{r}\}$. If $c\sim x_{1}$, then $G[X_{3}]$ has $x_{1}c$
as its only edge. Since $b\in\operatorname{EPN}(x_{1},X_{3})$ and $x_{2}%
\in\operatorname{EPN}(c,X_{3})$, $X_{3}$ is an $\operatorname{IR}(G)$-set
containing an edge, which is not the case. Hence $c\nsim x_{1}$, so $X_{3}$ is
an independent $\operatorname{IR}(G)$-set and $(X_{0}\overset{x_{1}b}{\sim
}_{H}X_{1}\overset{x_{2}c}{\sim}_{H}X_{2}\overset{bx_{1}}{\sim}_{H}%
X\overset{cx_{2}}{\sim}_{H}X_{0})$ is an induced $4$-cycle in~$H$%
.~$\blacksquare$

\bigskip

We use Corollary \ref{Cor_Not_independent} and Lemma \ref{Lem_Indep} to prove
our next result, which shows that not all connected graphs are
$\operatorname{IR}$-graphs.

\begin{proposition}
\label{Prop_diam2}If $H$ is an $\operatorname{IR}$-graph of diameter $2$, then
$H$ has an induced $C_{4}$.
\end{proposition}

\noindent\textbf{Proof.\hspace{0.1in}}Suppose, to the contrary, that $H$ is a
$C_{4}$-free graph with $\operatorname{diam}(H)=2$, but $H=G(\operatorname{IR}%
)$ for some graph $G$. By Corollary \ref{Cor_Not_independent}, each
$\operatorname{IR}$-set of $G$ is independent. Since $H$ is connected but not
complete, $\operatorname{IR}(G)\geq2$ and $G$ has at least three
$\operatorname{IR}$-sets. Say $G$ has $k$ $\operatorname{IR}$-sets. Let
$X_{0}=\{x_{1},...,x_{r}\}$ be any $\operatorname{IR}(G)$-set and, without
loss of generality, $X_{1}=\{b,x_{2},...,x_{r}\}$, where $b\sim x_{1}$, is an
$\operatorname{IR}(G)$-set such that $X_{0}\sim_{H}X_{1}$. As shown in the
second paragraph of the proof of Lemma \ref{Lem_Indep}, if $G$ has an
(independent) $\operatorname{IR}$-set $\{y,c,x_{3},...,x_{r}\}$, where
$y\in\{x_{1},b\}$, then $H$ contains an induced $C_{4}$. Hence all
$\operatorname{IR}(G)$-sets are independent sets of the form $Y_{i}%
=\{y_{i},x_{2},...,x_{r}\},\ i=1,...,k$, where, as shown in the first
paragraph of the proof of Lemma \ref{Lem_Indep}, $G[\{y_{1},...,y_{k}\}]$ is a
complete graph. But then $Y_{i}\sim_{H}Y_{j}$ for all $1\leq i<j\leq k$, which
implies that $H\cong G[\{y_{1},...,y_{k}\}]$ is also complete. This
contradicts $\operatorname{diam}(H)=2$.~$\blacksquare$

\bigskip

By Proposition~\ref{Prop_diam2}, trees of diameter $2$, i.e. the stars
$K_{1,n}$, are not $\operatorname{IR}$-graphs. We generalise this result to
include all non-complete graphs with universal vertices.

\begin{proposition}
\label{Prop_K1_n}If $H$ is a non-complete graph with a universal vertex, then
$H$ is not an $\operatorname{IR}$-graph.
\end{proposition}

\noindent\textbf{Proof.\hspace{0.1in}}Since $H$ is not complete and has a
universal vertex, $\operatorname{diam}(H)=2$. Suppose, contrary to the
statement, that $H=G(\operatorname{IR})$ for some graph $G$. Let $u$ be a
universal vertex of $H$ and $X=\{x_{1},...,x_{r}\}$ the $\operatorname{IR}%
(G)$-set corresponding to $u$. Since $\operatorname{diam}(G)=2$, Corollary
\ref{Cor_Number_IR-sets1} implies that $G[X]$ contains at most one edge. Since
$X$ does not lie on an induced $C_{4}$ in $H$, Corollary \ref{Cor_C4}$(i)$
implies that $X$ is independent.

Since $H$ is not complete, there exist $\operatorname{IR}(G)$-sets
$Y_{1},Y_{2},Y_{3}$ distinct from $X$ such that $Y_{1}\sim_{H}Y_{2}\sim
_{H}Y_{3}$ but $Y_{1}\nsim_{H}Y_{3}$. Moreover, $X\sim_{H}Y_{i}$ for
$i=1,2,3$. Without loss of generality, assume that $Y_{2}=\{b,x_{2}%
,...,x_{r}\}$, where $x_{1}\sim b$.

\begin{itemize}
\item We show that $Y_{1}=\{a,x_{2},...,x_{r}\}$, where $x_{1}\sim a$ and
$a\sim b$.
\end{itemize}

Since $Y_{1}\sim_{H}Y_{2}$, $|Y_{1}-Y_{2}|=1$. Suppose that $x_{1}\in Y_{1}$
and $b\notin Y_{1}$; say $Y_{1}=\{x_{1},y_{2},x_{3},...,x_{r}\}$, where
$y_{2}\neq b$ and $x_{2}\sim y_{2}$. Then $|Y_{1}-Y_{2}|=|\{x_{1},y_{2}\}|=2$,
a contradiction.

Suppose $\{x_{1},b\}\subseteq Y_{1}$. We may assume without loss of generality
that $Y_{1}=\{x_{1},b,x_{3},...,x_{r}\}$, where $x_{2}\sim b$. Then $Y_{1}$
and $Y_{2}$ differ only in that $x_{1}\in Y_{1}-Y_{2}$ and $x_{2}\in
Y_{2}-Y_{1}$. Since $Y_{1}\sim Y_{2}$, we deduce that $x_{1}\sim x_{2}$.
However, this is a contradiction because $X$ is independent.

Therefore $Y_{1}=\{a,x_{2},...,x_{r}\}$, where $x_{1}\sim a$ and $a\sim b$.
Similarly, $Y_{3}=\{c,x_{2},...,x_{r}\}$, where $c\neq a$, $x_{1}\sim c$ and
$c\sim b$. Since the symmetric difference $Y_{1}\vartriangle Y_{3}=\{a,c\}$
and $Y_{1}\nsim Y_{3}$, we deduce that $a\nsim c$. Since $|Y_{1}%
\cup\{c\}|>\operatorname{IR}(G)$, we know that $G[Y_{1}\cup\{c\}]$ contains an
edge. Since $ac\notin E(G)$ and $X$ is independent, we may assume without loss
of generality that $ax_{2}$ is an edge of $G[Y_{1}]$. (The argument is the
same if $c$ is adjacent to a vertex in $\{x_{2},...,x_{r}\}$.) But $Y_{1}$ is
an $\operatorname{IR}(G)$-set and $\operatorname{diam}(H)=2$, hence by
Corollary \ref{Cor_Number_IR-sets1}, $ax_{2}$ is the only edge of $G[Y_{1}]$.

Let $d\in\operatorname{EPN}(x_{2},Y_{1})$. Since $X$ is independent and $d\sim
x_{2}$, we know that $d\neq x_{1}$. Since $ax_{2}$ is the only edge of
$G[Y_{1}]$ and $d\in\operatorname{EPN}(x_{2},Y_{1})$, the set $Z=(Y_{1}%
-\{x_{2}\})\cup\{d\}=\{a,d,x_{3},...,x_{r}\}$ is independent, hence an
$\operatorname{IR}(G)$-set. However, $|X-Z|=2$, i.e., $d_{H}(X,Z)\geq2$, which
contradicts the fact that $X$ corresponds to the universal vertex $u$ of
$H$.~$\blacksquare$

\section{Trees with diameter $3$ or $4$}

\label{Sec_Trees}We continue our investigation of graphs that are realizable
as $\operatorname{IR}$-graphs by considering $\operatorname{IR}$-\emph{trees},
that is, trees that are $\operatorname{IR}$-graphs. Since all complete graphs
are $\operatorname{IR}$-graphs, $K_{1}$ and $K_{2}$ are $\operatorname{IR}%
$-trees. We will show (see Theorem \ref{Thm_Main}) that the smallest
non-complete $\operatorname{IR}$-tree is the double star $S(2,2)$, which has
order $6$. We know from Propositions \ref{Prop_diam2} and \ref{Prop_K1_n} that
trees of diameter $2$ (i.e.~stars) are not $\operatorname{IR}$-graphs. We now
focus on trees of diameter $3$ and $4$. Lemma \ref{Lem_Not_Path} below,
another result that explores the role of external private neighbours and
flip-sets, is useful in both cases as it gives information on the structure of
$\operatorname{IR}$-graphs of graphs with certain types of $\operatorname{IR}$-sets.

\begin{lemma}
\label{Lem_Not_Path}Let $H$ be a connected $\operatorname{IR}$-graph of a
graph $G$. Suppose $X$ is an $\operatorname{IR}(G)$-set such that exactly
three vertices $x_{1},x_{2},x_{3}$ have positive degree in $G[X]$. For
$i=1,2,3$, let $x_{i}^{\prime}\in\operatorname{EPN}(x_{i},X)$ and let
$X^{\prime}$ be the flip-set of $X$ using $\{x_{1}^{\prime},x_{2}^{\prime
},x_{3}^{\prime}\}$. If $d_{H}(X,X^{\prime})=3$, then $H$ contains an induced
$4$-cycle, or an induced double star $S(2,2)$ of which $X$ and $X^{\prime}$
are antipodal leaves.
\end{lemma}

\noindent\textbf{Proof.}\hspace{0.1in}Assume $d_{H}(X,X^{\prime})=3$ and $H$
does not contain an induced $C_{4}$. Let $P:(X=X_{0},X_{1},X_{2}%
,X_{3}=X^{\prime})$ be an $X$-$X^{\prime}$ geodesic in $H$. Since
$|X-X^{\prime}|=d_{H}(X,X^{\prime})=3$, we may assume without loss of
generality that $X_{1}=\{x_{1}^{\prime},x_{2},...,x_{r}\}$, $X_{2}%
=\{x_{1}^{\prime},x_{2}^{\prime},x_{3},...,x_{r}\}$ and $X_{3}=\{x_{1}%
^{\prime},x_{2}^{\prime},x_{3}^{\prime},x_{4},...,x_{r}\}$. Since
$x_{1}^{\prime}\nsim x_{2},x_{3},...,x_{r}$ by the private neighbour property
of $x_{1}$, the only possible edge in $G[X_{1}]$ is $x_{2}x_{3}$. By Corollary
\ref{Cor_C4}$(i)$ and our assumption that $H$ does not contain an induced
$C_{4}$, we may assume $x_{2}\nsim x_{3}$. Since $x_{1},x_{2},x_{3}$ have
positive degree in $G[X]$, $x_{1}\sim x_{2},x_{3}$. Similarly, the only
possible edge in $G[X_{2}]$ is $x_{1}^{\prime}x_{2}^{\prime}$, and again we
may assume that $x_{1}^{\prime}\nsim x_{2}^{\prime}$. Also applying Corollary
\ref{Cor_C4}$(i)$ to $X_{3}$, we have either (a) $x_{1}^{\prime}\nsim
x_{3}^{\prime}\nsim x_{2}^{\prime}$ or (b) $x_{1}^{\prime}\sim x_{3}^{\prime
}\sim x_{2}^{\prime}$.

Suppose (a) holds and consider $R=\{x_{1},x_{2},x_{3}^{\prime},x_{4}%
,...,x_{r}\}$. Since $x_{3}^{\prime}\in\operatorname{EPN}(x_{3},X)$,
$x_{1}x_{2}$ is the only edge of $G[R]$. However, $x_{i}^{\prime}%
\in\operatorname{EPN}(x_{i},R)$ for $i=1,2$, so $R$ is an $\operatorname{IR}%
(G)$-set. By Corollary \ref{Cor_C4}$(i)$, $H$ has a $4$-cycle, which is not
the case.

Hence assume (b) holds and consider the set $W=\{x_{2}^{\prime},x_{3}^{\prime
},x_{3},x_{4},...,x_{r}\}$.

\begin{itemize}
\item[$\ast$] Since $x_{3}\sim x_{1}$ while $x_{1}$ is nonadjacent to
$x_{2}^{\prime},x_{3}^{\prime},x_{4},...,x_{r}$,\hspace{0.1in}$x_{1}%
\in\operatorname{EPN}(x_{3},W)$.

\item[$\ast$] Since $x_{3}^{\prime}\sim x_{1}^{\prime}$ while $x_{1}^{\prime}$
is nonadjacent to $x_{2}^{\prime},x_{3},...,x_{r}$,\hspace{0.1in}%
$x_{1}^{\prime}\in\operatorname{EPN}(x_{3}^{\prime},W)$.

\item[$\ast$] Since $x_{2}\sim x_{2}^{\prime}$ while $x_{2}$ is nonadjacent to
$x_{3}^{\prime},x_{3},...,x_{r}$,\hspace{0.1in}$x_{2}\in\operatorname{EPN}%
(x_{2}^{\prime},W)$.

\item[$\ast$] Since $x_{i}$ is isolated in $G[W]$ for all $i\geq4$, $x_{i}%
\in\operatorname{PN}(x_{i},W)$.
\end{itemize}

\noindent Therefore $W$ is an $\operatorname{IR}(G)$-set. Since $x_{3}%
^{\prime}\sim x_{1}^{\prime}$, $W\sim_{H}X_{2}$. By Proposition
\ref{Prop_IR_sets}, the flip-set $U=\{x_{1},x_{1}^{\prime},x_{2}%
,x_{4},...,x_{r}\}$ of $W$ is also an $\operatorname{IR}(G)$-set. Since
$x_{1}\sim x_{3}$, $U\sim_{H}X_{1}$. Since $U$ and $W$ are nonadjacent to each
other and to the $\operatorname{IR}(G)$-sets $X,X_{2},X_{3}$ and
$X,X_{1},X_{3}$, respectively, $H$ contains the double star $S(2,2)$ as
induced subgraph. Since $P$ is an $X$-$X^{\prime}$ geodesic, $X$ and
$X^{\prime}$ are antipodal leaves of this double star.~$\blacksquare$

\bigskip

We continue our investigation of trees that are realizable as
$\operatorname{IR}$-trees and determine the smallest $\operatorname{IR}$-tree
with diameter $3$. As stated in the introduction, a frequently studied
property of a reconfiguration graph is its diameter. In particular, a small
diameter of an $\operatorname{IR}$-graph indicates that it is relatively easy
to transition from one $\operatorname{IR}$-set in the source graph to any
other one. \label{Wednesday} We conjecture (see Conjecture
\ref{Con_doublestar}) that the double stars $S(2k,2k),\ k\geq1$, are the only
$\operatorname{IR}$-trees with diameter~$3$.

\begin{proposition}
\label{Prop_Diam3}The double star $S(2,2)$ is the unique smallest
$\operatorname{IR}$-tree with diameter~$3$.
\end{proposition}

\begin{figure}[ptb]%
\centering
\includegraphics[
height=1.4935in,
width=4.2004in
]%
{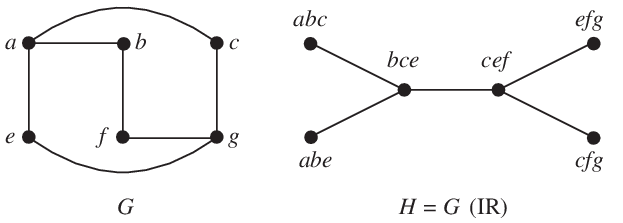}%
\caption{A graph $G$ and its $\operatorname{IR}$-graph $H=S(2,2)$}%
\label{Fig_S22}%
\end{figure}

\noindent\textbf{Proof.\hspace{0.1in}}Suppose $T$ with $\operatorname{diam}%
(T)=3$ is an $\operatorname{IR}$-tree of a graph $G$. By Corollaries
\ref{Cor_C4} and \ref{Cor_Number_IR-sets1}, all $\operatorname{IR}$-sets of
$G$ are independent or induce a graph that has exactly three vertices of
positive degree. If all $\operatorname{IR}(G)$-sets are independent, then the
$\operatorname{IR}$-graph of $G$ has a cycle, by Lemma \ref{Lem_Indep}. On the
other hand, if $G$ has an $\operatorname{IR}$-set $X$ containing exactly three
vertices $x_{1},x_{2},x_{3}$ of positive degree in $G[X]$, then Lemma
\ref{Lem_Not_Path} implies that $T$ has $S(2,2)$ as subgraph. The graph $G$ in
Figure~\ref{Fig_S22} is an example of a graph for which $G(\operatorname{IR}%
)\cong S(2,2)$. (Verifying this involves an exhaustive but straightforward
search for $\operatorname{IR}(G)$-sets and their adjacencies.)~$\blacksquare$

\bigskip%

%TCIMACRO{\FRAME{ftbpFU}{4.2004in}{1.4935in}{0pt}{\Qcb{A graph $G$ and its
%$\operatorname{IR}$-graph $H=S(2,2)$}}{\Qlb{Fig_S22}}{ir_s22.eps}%
%{\special{ language "Scientific Word";  type "GRAPHIC";
%maintain-aspect-ratio TRUE;  display "USEDEF";  valid_file "F";
%width 4.2004in;  height 1.4935in;  depth 0pt;  original-width 4.1502in;
%original-height 1.4581in;  cropleft "0";  croptop "1";  cropright "1";
%cropbottom "0";  filename 'IR_S22.eps';file-properties "XNPEU";}} }%
%BeginExpansion
%EndExpansion

For the remainder of the paper we consider trees of diameter $4$. Our aim is
to show that the double spider $\mathrm{Sp}(1,1;1,2)$ is the unique smallest
$\operatorname{IR}$-tree with diameter $4$. This case turns out to be more
challenging than the previous topics and requires a technical lemma. We first
state a simple observation for referencing.

\begin{observation}
\label{Ob_diam5}If $P:(v_{0},...,v_{4})$ is a path in a tree $T$ and $u$ is a
vertex such that $d(u,v_{2})\geq3$ or $d(u,v_{1})\geq4$, then
$\operatorname{diam}(T)\geq5$.
\end{observation}

Lemma \ref{Lem_Not_Path} concerns an $\operatorname{IR}(G)$-set $X$ such that
exactly three vertices $x_{1},x_{2},x_{3}\in X$ have positive degree in $G[X]$
and such that the distance, in $G(\operatorname{IR})$, between $X$ and its
flip-set $X^{\prime}$ is $3$. We now consider the same situation except that
$d_{G(\operatorname{IR})}(X,X^{\prime})=4$; this apparently small change
results in a much more complex situation. Lemma \ref{Lem_Tec} below allows us
to characterize the unique smallest $\operatorname{IR}$-tree with diameter
$4$. 

\begin{lemma}
\label{Lem_Tec}Let $H$ be an $\operatorname{IR}$-graph of $G$. Suppose $X$ is
an $\operatorname{IR}(G)$-set such that exactly three vertices $x_{1}%
,x_{2},x_{3}$ have positive degree in $G[X]$. For $i=1,2,3$, let
$x_{i}^{\prime}\in\operatorname{EPN}(x_{i},X)$ and let $X^{\prime}$ be the
flip-set of $X$ using $\{x_{1}^{\prime},x_{2}^{\prime},x_{3}^{\prime}\}$.
Suppose $d_{H}(X,X^{\prime})=4$ and $P:(X=X_{0},...,X_{4}=X^{\prime})$ is an
$X-X^{\prime}$ geodesic. Then

\begin{enumerate}
\item[$(i)$] if $x_{j}^{\prime}\in X_{i}$ for some $i$, then $x_{j}^{\prime
}\in X_{\ell}$ for all $\ell=i,...,4$;

\item[$(ii)$] there exists exactly one vertex $a\in\bigcup_{i=1}^{3}X_{i}$
such that $a\notin X\cup\{x_{1}^{\prime},x_{2}^{\prime},x_{3}^{\prime}\}$;

\item[$(iii)$] for this vertex $a$, if $X_{i}\overset{x_{j}a}{\sim}X_{i+1}$
for some $i$ and $j$, then $X_{\ell}\overset{ax_{j}^{\prime}}{\sim}X_{\ell+1}$
for some $\ell\geq i+1$, and $G[\{x_{j},a,x_{j}^{\prime}\}]=K_{3}$;

\item[$(iv)$] $\{x_{4},...,x_{r}\}$ is a subset of $X_{i}$ for each
$i=0,...,4$;

\item[$(v)$] if $H$ is a tree such that $|V(H)|\leq7$ and $\operatorname{diam}%
(H)=4$, then at most one $x_{i},\ i=1,2,3$, has an $X$-external private
neighbour $y_{i}\neq x_{i}^{\prime}$; if $y_{i}$ exists, then $y_{i}=a$,
$a\sim x_{i}^{\prime}$ and the flip-set $(X^{\prime}-\{x_{i}^{\prime}%
\})\cup\{a\}$ of $X$ using $a$ instead of $x_{i}^{\prime}$ is $X_{3}$.
\end{enumerate}
\end{lemma}

\noindent\textbf{Proof.\hspace{0.1in}}To obtain $X^{\prime}$ from $X$ requires
exactly four swaps, and in three of these the $x_{i}^{\prime},\ i=1,2,3$, are
swapped in, while some vertex $a\notin\{x_{1}^{\prime},x_{2}^{\prime}%
,x_{3}^{\prime}\}$ is swapped in during another swap (but not the last,
obviously). Thus, if $x_{i}^{\prime}$ has been swapped into $X_{j}$ for some
$j\geq1$, then $x_{i}^{\prime}$ is never swapped out, otherwise it would have
to be swapped in again, necessitating too many swaps. If some $x_{j}$, $j>3$,
is swapped for a vertex $v$, then $v\notin\{x_{1},x_{2},x_{3},x_{1}^{\prime
},x_{2}^{\prime},x_{3}^{\prime}\}$ (since $x_{j}$, being isolated in $G[X]$
and by the private neighbour property, is nonadjacent to all of these
vertices), so $x_{j}$ has to be swapped in again, again resulting in too many
swaps. This proves $(i)$, $(ii)$ and$~(iv)$.

Suppose vertex $x_{j}$ is swapped out for $a$. To avoid having too many swaps,
vertex $a$ is later swapped out for some $x_{i}^{\prime}$. If $i\neq j$, then
$x_{\ell}$ is swapped out for $x_{j}^{\prime}$, where $\ell\neq j$. But since
$x_{j}^{\prime}\in\operatorname{EPN}(x_{j},X)$, $x_{\ell}\nsim x_{j}^{\prime}%
$, a contradiction. Hence $i=j$ and $x_{j}\sim_{G}a\sim_{G}x_{j}^{\prime}%
\sim_{G}x_{j}$. This proves~$(iii)$.

\noindent$(v)\hspace{0.1in}$Assume $H$ is a tree such that $|V(H)|\leq7$ and
$\operatorname{diam}(H)=4$. Suppose, for some $i=1,2,3$, $x_{i}$ has an
$X$-external private neighbour $y_{i}\neq x_{i}^{\prime}$. Let $Y_{i}$ be the
flip-set of $X$ using $y_{i}$ instead of $x_{i}^{\prime}$, that is,
$Y_{i}=\{y_{i},x_{j}^{\prime},x_{k}^{\prime},x_{4},...,x_{r}\}$,
$\{i,j,k\}=\{1,2,3\}$. By Proposition \ref{Prop_IR_sets}, $Y_{i}$ is an
$\operatorname{IR}(G)$-set.

\begin{itemize}
\item Suppose $y_{i}\nsim x_{i}^{\prime}$. By the second condition in $(iii)$,
$y_{i}\neq a$, and since $y_{i}\notin X\cup\{x_{1}^{\prime},x_{2}^{\prime
},x_{3}^{\prime}\}$, $Y_{i}\notin\{X_{0},...,X_{4}\}$. But now $\{x_{j}%
,y_{i},x_{i}^{\prime},x_{4},...,x_{r}\}$ and $\{x_{k},y_{i},x_{i}^{\prime
},x_{4},...,x_{r}\}$ are independent $\operatorname{IR}(G)$-sets different
from $Y_{i},X_{0},...,X_{4}$ and $T$ has order at least eight, a
contradiction. Therefore $y_{i}\sim x_{i}^{\prime}$, so $Y_{i}\sim_{H}X_{4}$.

\item If $Y_{i}\neq X_{3}$, then either $d_{H}(X,Y_{i})=5$ (if $Y_{i}$ is
nonadjacent to $X_{0},...,X_{3}$) or $H$ has a cycle (otherwise), which is
impossible. Hence $Y_{i}=X_{3}$.

\item Since $y_{i}\notin X\cup\{x_{1}^{\prime},x_{2}^{\prime},x_{3}^{\prime
}\}$, $(ii)$ implies that $y_{i}=a$.

\item We still have to show that no $x_{j}$, where $j\neq i$, has an
$X$-external private neighbour. Suppose $y_{j}\in\operatorname{EPN}%
(x_{j},X)-\{x_{j}^{\prime}\},\ j\neq i$. As shown above, $y_{j}=a$. Now we
have $a\in\operatorname{EPN}(x_{i},X)\cap\operatorname{EPN}(x_{j},X)$, a
contradiction.~$\blacksquare$
\end{itemize}

\label{P5}%
%TCIMACRO{\FRAME{ftbpFU}{4.4209in}{1.8628in}{0pt}{\Qcb{A graph $F$ and its
%$\operatorname{IR}$-graph, the double spider $\QTR{rm}{Sp}(1,1;1,2)$%
%}}{\Qlb{Fig_Tree}}{ir_conn.eps}{\special{ language "Scientific Word";
%type "GRAPHIC";  maintain-aspect-ratio TRUE;  display "USEDEF";
%valid_file "F";  width 4.4209in;  height 1.8628in;  depth 0pt;
%original-width 4.6241in;  original-height 1.932in;  cropleft "0";
%croptop "1";  cropright "1";  cropbottom "0";
%filename 'IR_Conn.eps';file-properties "XNPEU";}} }%
%BeginExpansion
\begin{figure}[ptb]%
\centering
\includegraphics[
height=1.8628in,
width=4.4209in
]%
{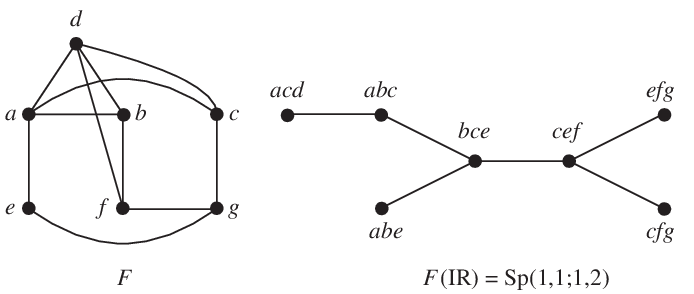}%
\caption{A graph $F$ and its $\operatorname{IR}$-graph, the double spider
$\mathrm{Sp}(1,1;1,2)$}%
\label{Fig_Tree}%
\end{figure}
%EndExpansion

Concerning the realizability of trees as $\operatorname{IR}$-trees, we next
consider trees with diameter $4$ and show that $\mathrm{Sp}(1,1;1,2)$ is the
unique smallest $\operatorname{IR}$-tree with diameter $4$. The proof has two
main cases, depending on whether the source graph $G$ has an
$\operatorname{IR}$-set $X$ such that $G[X]$ has four non-isolated vertices or not.

\begin{theorem}
\label{Thm_P5}The double spider $\mathrm{Sp}(1,1;1,2)$ is the unique smallest
$\operatorname{IR}$-tree with diameter $4$.
\end{theorem}

\noindent\textbf{Proof.\hspace{0.1in}}As illustrated in Figure \ref{Fig_Tree},
the double spider $\mathrm{Sp}(1,1;1,2)$ is an $\operatorname{IR}$-tree with
diameter $4$ and order $7$. (Again, verifying this involves an exhaustive but
straightforward search for $\operatorname{IR}(G)$-sets. The graph $F$ is
obtained from the graph $G$ in Figure \ref{Fig_S22} by adding a new vertex
$d$, joining it to $a,b,c$ and $f$.)

Let $T$ be a tree with $\operatorname{diam}(T)=4$ and $|V(T)|\leq7$, and
suppose $T$ is the $\operatorname{IR}$-graph of a graph $G$. We show that $T$
contains $S(2,2)$. Since $\operatorname{diam}(T)=4$ and $\operatorname{diam}%
(S(2,2))=3$, $T$ has order at least $7$, and the only possibility is $T\cong
S(1,1;1,2\}$.

Let $P:(X_{0},...,X_{4})$ be a diametrical path in $T$. By Corollary
\ref{Cor_Number_IR-sets1}, the $\operatorname{IR}$-sets of $G$ induce graphs
that have at most four vertices of positive degree. Let $X=\{x_{1}%
,...,x_{r}\}$ be an $\operatorname{IR}(G)$-set that induces a subgraph with
the largest number of vertices of positive degree. Say $x_{1},...,x_{k}$,
$k\leq4$, have positive degree in $G[X]$. Let $x_{i}^{\prime}\in
\operatorname{EPN}(x_{i},X)$ and $X^{\prime}=(X-\{x_{1},...,x_{k}%
\})\cup\{x_{1}^{\prime},...,x_{k}^{\prime}\}$; note that $|X-X^{\prime}|=k$.
By Proposition \ref{Prop_IR_sets}, $X^{\prime}$ is an $\operatorname{IR}%
(G)$-set. By Corollary \ref{Cor_C4}$(i)$ and Lemma \ref{Lem_Indep},
$k\in\{3,4\}$. We consider the two cases separately.

\noindent\textbf{Case 1\hspace{0.1in}}$k=4$. Then $d_{T}(X,X^{\prime})=4$,
hence without loss of generality,
\[
X_{0}=X,X_{1}=\{x_{1}^{\prime},x_{2},...,x_{r}\},X_{2}=\{x_{1}^{\prime}%
,x_{2}^{\prime},x_{3},...,x_{r}\},X_{3}=\{x_{1}^{\prime},x_{2}^{\prime}%
,x_{3}^{\prime},x_{4},...,x_{r}\}\ \text{and}\ X_{4}=X^{\prime}.
\]

Since $x_{1}^{\prime}\in\operatorname{EPN}(x_{1},X)$, $x_{1}^{\prime}$ is
nonadjacent to all of $x_{2},...,x_{r}$, hence $x_{2},x_{3},x_{4}$ are the
only possible vertices with positive degree in $G[X_{1}]$. By Corollary
\ref{Cor_C4}, either all or none of them have positive degree. We consider the
two subcases separately.\medskip

\noindent\textbf{Case 1.1\hspace{0.1in}}Say $x_{2},x_{3},x_{4}$ all have
positive degree in $G[X_{1}]$. Let $x_{i}^{\prime\prime}\in\operatorname{EPN}%
(x_{i},X_{1}),\ i=2,3,4$, and let $X^{\prime\prime}$ be the flip-set of
$X_{1}$ using $\{x_{2}^{\prime\prime},x_{3}^{\prime\prime},x_{4}^{\prime
\prime}\}$. Then $d_{T}(X_{1},X^{\prime\prime})\geq3$. If $d_{T}%
(X_{1},X^{\prime\prime})\geq4$, then $\operatorname{diam}(T)>4$ by Observation
\ref{Ob_diam5}, a contradiction. Hence $d_{T}(X_{1},X^{\prime\prime})=3$. By
Lemma \ref{Lem_Not_Path}, $T$ contains $S(2,2)$. Since $\operatorname{diam}%
(T)=4$ and $T$ has order $7$, whereas $S(2,2)$ has diameter $3$ and order $6$,
$T\cong S(1,1;1,2\}$.\medskip

\noindent\textbf{Case 1.2\hspace{0.1in}}$\{x_{2},x_{3},x_{4}\}$ is
independent. Since $\{v\in X:\deg_{G[X]}(v)>0\}=\{x_{1},...,x_{4}\}$, we
deduce that $x_{1}$ is adjacent to each of $x_{2},x_{3}$ and $x_{4}$.
Moreover, since $x_{1}^{\prime}\in\operatorname{EPN}(x_{1},X)$, the set
$X_{1}$ is independent. Again by the private neighbour property, the only
possible edge in $G[X_{2}]$ is $x_{1}^{\prime}x_{2}^{\prime}$. But by
Corollary \ref{Cor_C4}$(i)$ applied to $X_{2}$, $x_{1}^{\prime}\nsim
x_{2}^{\prime}$, hence $X_{2}$ is independent. By Corollary \ref{Cor_C4}$(ii)$
applied to $X_{1}$, at most one vertex in $X_{1}$ has an external private
neighbour. Since $x_{2}^{\prime}\in\operatorname{EPN}(x_{2},X)$ and
$x_{1}^{\prime}\nsim x_{2}^{\prime}$, $x_{2}^{\prime}\in\operatorname{EPN}%
(x_{2},X_{1})$. Hence $x_{2}$ is the unique vertex in $X_{1}$ with external
private neighbours. But for $i=3,4$, $x_{i}^{\prime}\in\operatorname{EPN}%
(x_{i},X)$. We deduce that $x_{1}^{\prime}\sim x_{3}^{\prime},x_{4}^{\prime}$
to ensure that $\operatorname{EPN}(x_{3},X_{1})=\varnothing=\operatorname{EPN}%
(x_{4},X_{1})$. By Corollary \ref{Cor_C4}$(i)$ applied to $X_{3}$, in which
$x_{4}$ is isolated (by the private neighbourhood property and because
$x_{5},...,x_{r}$ are isolated in $G[X]$), $x_{2}^{\prime}\sim x_{3}^{\prime}$.

Note that $x_{i}\in\operatorname{EPN}(x_{i}^{\prime},X_{3})$ for $i=2,3$, but
$x_{1}\sim x_{4}$, so $x_{1}\notin\operatorname{EPN}(x_{1}^{\prime},X_{3})$.
Let $u\in\operatorname{EPN}(x_{1}^{\prime},X_{3})$ ($u$ exists because
$x_{1}^{\prime}$ is not isolated in $G[X_{3}]$) and let $B$ be the flip-set of
$X_{3}$ using $\{u,x_{2},x_{3}\}$. Then $B=\{u,x_{2},...,x_{r}\}$ is an
$\operatorname{IR}(G)$-set, by Proposition \ref{Prop_IR_sets}. Since $u\sim
x_{1}^{\prime}$, $B\sim_{T}X_{1}$. Since $T$ is a tree and $B\neq X_{2}$,
$d_{T}(B,X_{3})=3$. By Lemma \ref{Lem_Not_Path}, $T$ contains $S(2,2)$, and as
in Case 1.1, $T\cong S(1,1;1,2\}$.\medskip

\noindent\textbf{Case 2\hspace{0.1in}}$k=3$, that is, $x_{1},x_{2},x_{3}$ are
the only vertices in $X$ that have positive degree in $G[X]$, and
$|X-X^{\prime}|=3$. Therefore $3\leq d_{T}(X,X^{\prime})\leq4$. If
$d_{T}(X,X^{\prime})=3$, then Lemma \ref{Lem_Not_Path} implies that $T$
contains $S(2,2)$ and we are done. Hence we assume that $d_{T}(X,X^{\prime
})=4$ and that $P:(X=X_{0},...,X_{4}=X^{\prime})$ is an $X-X^{\prime}$
geodesic. We may also assume without loss of generality (otherwise we can just
relabel) that $x_{1}^{\prime}$ is swapped into the irredundant set before
$x_{2}^{\prime}$, which, in turn, is swapped in before $x_{3}^{\prime}$. By
Lemma \ref{Lem_Tec}$(ii)$, $\bigcup_{i=0}^{4}X_{i}$ contains exactly one
vertex $a\notin X\cup\{x_{1}^{\prime},x_{2}^{\prime},x_{3}^{\prime}\}$.

To obtain $X^{\prime}$ from $X$ requires exactly four steps, and, as shown in
Lemma \ref{Lem_Tec}, in three of these the $x_{i}^{\prime},\ i=1,2,3$, are
swapped in, while $a$ is swapped in during another step (but not the last,
obviously). We consider the possibilities for the step in which $a$ is
swapped.\medskip

\noindent\textbf{Case 2.1\hspace{0.1in}}Vertex $a$ is swapped in first. Then
step 1 is either $X_{0}\overset{x_{1}a}{\sim}X_{1}$, $X_{0}\overset{x_{2}%
a}{\sim}X_{1}$ or $X_{0}\overset{x_{3}a}{\sim}X_{1}$.

\noindent\underline{Suppose step 1 is $X_{0}\overset{x_{1}a}{\sim}X_{1}$.}

By our assumption above on the order in which the $x_{i}^{\prime}$ are
swapped, step 2 is $X_{1}\overset{ax_{1}^{\prime}}{\sim}X_{2}$. By Lemma
\ref{Lem_Tec}$(iii)$, $G[\{x_{1},a,x_{1}^{\prime}\}]=K_{3}$, which implies
that $X_{0}\sim_{T}X_{1}\sim_{T}X_{2}\sim_{T}X_{0}$ and $T$ has a cycle, which
is not the case.

\noindent\underline{Suppose step 1 is $X_{0}\overset{x_{2}a}{\sim}X_{1}$.}

Then step 2 is $X_{1}\overset{x_{1}x_{1}^{\prime}}{\sim}X_{2}$ and step 3 is
$X_{2}\overset{ax_{2}^{\prime}}{\sim}X_{3}$. Hence
\[
X_{1}=\{x_{1},a,x_{3},...,x_{r}\},\ X_{2}=\{x_{1}^{\prime},a,x_{3}%
,...,x_{r}\},\ X_{3}=\{x_{1}^{\prime},x_{2}^{\prime},x_{3},...,x_{r}%
\},\ X_{4}=\{x_{1}^{\prime},x_{2}^{\prime},x_{3}^{\prime},...,x_{r}\}.
\]
We first show that

\begin{enumerate}
\item[(a)] $a\notin\operatorname{EPN}(x_{2},X)$, hence $a$ is not isolated in
$G[X_{1}]$, and

\item[(b)] $a\ $is\ adjacent\ to$\ x_{1}$\ or$\ x_{3},\ $but\ not\ to$\ x_{4}%
,...,x_{r}.$
\end{enumerate}

\noindent We then deduce that

\begin{enumerate}
\item[(c)] $T$ contains $S(2,2).$
\end{enumerate}

\noindent(a)\hspace{0.1in}If $a\in\operatorname{EPN}(x_{2},X)$, then the
flip-set of $X$ using $a$ instead of $x_{2}^{\prime}$, which is $\{x_{1}%
^{\prime},a,x_{3}^{\prime},...,x_{r}\}$, equals $X_{3}$ by Lemma
\ref{Lem_Tec}$(v)$. However, $a\notin X_{3}$ and we have a contradiction.
Therefore (a) holds. This implies that $a$ is not isolated in $G[X_{1}%
]$.\smallskip

\noindent(b)\hspace{0.1in}Suppose $a$ is adjacent to one of $x_{4},...,x_{r}$,
say $x_{4}$. This means that $a$ and $x_{4}$ have nonempty $X_{1}$-external
private neighbourhoods. Then $x_{1}\nsim x_{3}$, otherwise $G[X_{1}]$ has four
vertices of positive degree, which is not the case since $k=3$.%
%TCIMACRO{\FRAME{ftbpFU}{6.5535in}{1.4762in}{0pt}{\Qcb{Case 2.1(b): The
%$\operatorname{IR}$-sets $X,\ X_{1}$ and $Q_{1}$ when $a$ is adjacent to
%$x_{2},\ x_{4}$ and $x_{5}$}}{\Qlb{Fig_551}}{thm_5_5_1.eps}%
%{\special{ language "Scientific Word";  type "GRAPHIC";
%maintain-aspect-ratio TRUE;  display "USEDEF";  valid_file "F";
%width 6.5535in;  height 1.4762in;  depth 0pt;  original-width 7.2506in;
%original-height 1.6094in;  cropleft "0";  croptop "1";  cropright "1";
%cropbottom "0";  filename 'Thm_5_5_1.eps';file-properties "XNPEU";}} }%
%BeginExpansion
\begin{figure}[ptb]%
\centering
\includegraphics[
height=1.4762in,
width=6.5535in
]%
{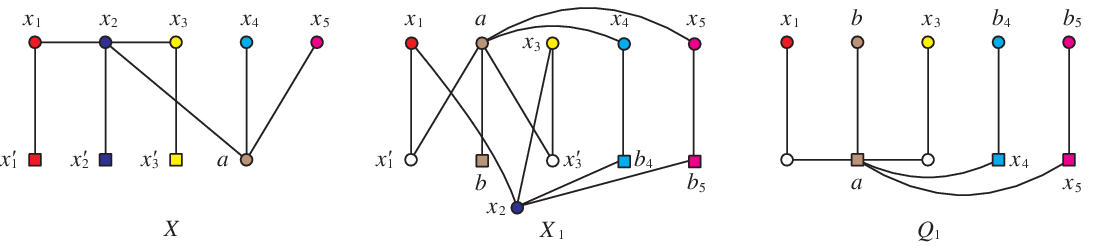}%
\caption{Case 2.1(b): The $\operatorname{IR}$-sets $X,\ X_{1}$ and $Q_{1}$
when $a$ is adjacent to $x_{2},\ x_{4}$ and $x_{5}$}%
\label{Fig_551}%
\end{figure}
%EndExpansion

$\bullet\hspace{0.1in}$Assume first that $a\nsim x_{1},x_{3}$. Then $x_{1}$
and $x_{3}$ are isolated in $G[X_{1}]$. By Corollary \ref{Cor_C4}$(i)$,
$ax_{4}$ is not the only edge of $G[X_{1}]$. Suppose $a\sim_{G}x_{5}$. (See
Figure \ref{Fig_551}.) Then $a,x_{4},x_{5}$ have positive degree in $G[X_{1}%
]$. Let $b,b_{4},b_{5}$ be $X_{1}$-external private neighbours of
$a,x_{4},x_{5}$, respectively and let $Q_{1}$ be the flip-set of $X_{1}$ using
$\{b_{1},b_{4},b_{5}\}$. Then $Q_{1}=\{x_{1},x_{3},b,b_{4},b_{5},...,x_{r}\}$
is an $\operatorname{IR}(G)$-set (Proposition \ref{Prop_IR_sets}). But
$\{x_{1},b_{4},b_{5}\}\subseteq Q_{1}-X_{2}$, hence $d_{T}(Q_{1},X_{2})\geq3$
and $\operatorname{diam}(T)\geq5$, by Observation \ref{Ob_diam5}. Therefore
$a\nsim x_{5},...,x_{r}$.

The only possibility therefore is either $a\sim x_{1}$ or $a\sim x_{3}$ (not
both, since $k=3$).

$\bullet\hspace{0.1in}$Say $a\sim x_{1}$. Then $x_{3}$ is isolated in
$G[X_{1}]$. In $G[X_{2}]$, $x_{1}^{\prime}\nsim x_{3},...,x_{r}$ and $a\nsim
x_{3},x_{5},...,x_{r}$. Using Corollary \ref{Cor_C4}$(i)$ and the fact that
$a\sim x_{4}$, we see that $a\sim x_{1}^{\prime}$. This means that
$x_{1}^{\prime}\notin\operatorname{EPN}(x_{1},X_{1})$. Let $c_{1},c,c_{4}$ be
$X_{1}$-external private neighbours of $x_{1},a,x_{4}$, respectively and let
$Q_{2}$ be the flip-set of $X_{1}$ using $\{c_{1},c,c_{4}\}$. Then
$Q_{2}=\{c_{1},c,x_{3},c_{4},x_{5},...,x_{r}\}$. Since $c_{1}\neq
x_{1}^{\prime}$, $\{c_{1},c,c_{4}\}\subseteq Q_{2}-X_{2}$, hence $d_{T}%
(Q_{2},X_{2})\geq3$ and $\operatorname{diam}(T)\geq5$, by Observation
\ref{Ob_diam5}.

$\bullet\hspace{0.1in}$Say $a\sim x_{3}$. Then $x_{1}$ is isolated in
$G[X_{1}]$. Let $d,d_{3},d_{4}$ be $X_{1}$-external private neighbours of
$a,x_{3},x_{4}$, respectively and let $Q_{3}$ be the flip-set of $X_{1}$ using
$\{d,d_{3},d_{4}\}$. Then $Q_{3}=\{x_{1},d,d_{3},d_{4},x_{5},...,x_{r}\}$, and
a contradiction follows as above.

We have now proved that $a$ is not adjacent to $x_{4},...,x_{r}$. Since $a$ is
not isolated in $G[X_{1}]$, (b) follows.\smallskip

\noindent(c)\hspace{0.1in}Corollary \ref{Cor_C4}$(i)$ and (b), together with
the assumptions of Case 2, ensure that $x_{1},x_{3}$ and $a$ have positive
degrees in $G[X_{1}]$. Let $y_{1},y_{2},y_{3}$ be $X_{1}$-external private
neighbours of $x_{1},a,x_{3}$, respectively, and let $Q_{4}=\{y_{1}%
,y_{2},y_{3},x_{4},...,x_{r}\}$ be the resulting flip-set of $X_{1}$. Since
$\operatorname{diam}(T)=4$ and $|Q_{4}-X_{1}|=3$, we deduce that $3\leq
d_{T}(X_{1},Q_{4})\leq4$. But if $d_{T}(X_{1},Q_{4})=4$, then, by Observation
\ref{Ob_diam5}, $\operatorname{diam}(T)\geq5$. Therefore $d_{T}(X_{1}%
,Q_{4})=3$. Applying Lemma \ref{Lem_Not_Path} to $X_{1}$, we deduce that $T$
contains $S(2,2)$.\medskip\ 

\noindent\underline{Suppose step 1 is $X_{0}\overset{x_{3}a}{\sim}X_{1}$.}

Then%
\[
X_{1}=\{x_{1},x_{2},a,x_{4},...,x_{r}\},X_{2}=\{x_{1}^{\prime},x_{2}%
,a,x_{4},...,x_{r}\}\ \text{and\ }X_{3}=\{x_{1}^{\prime},x_{2}^{\prime
},a,x_{4},...,x_{r}\}.
\]

$\bullet\hspace{0.1in}$Suppose $\,a\in\operatorname{EPN}(x_{3},X)$. We show that

\begin{enumerate}
\item[(d)] $x_{1}^{\prime},x_{2}^{\prime},a$ have positive degree in $X_{3}$,
and $X_{3}$-private neighbours $x_{1},x_{2},x_{3}$, respectively.
\end{enumerate}

Since $a\in\operatorname{EPN}(x_{3},X)$, $a$ is isolated in $G[X_{1}]$. By
Corollary \ref{Cor_C4}$(i)$, $x_{1}\nsim x_{2}$, that is, $X_{1}$ is
independent. Now the only possible edge in $G[X_{2}]$ is $ax_{1}^{\prime}$,
and again we deduce that $X_{2}$ is independent. If $a\nsim x_{2}^{\prime}$,
then $Q_{5}=\{x_{1},x_{2}^{\prime},a,x_{4},...,x_{r}\}$ is an independent
$\operatorname{IR}(G)$-set different from $X_{2}$ such that $X_{1}%
\overset{x_{2}x_{2}^{\prime}}{\sim}Q_{5}\overset{x_{1}x_{1}^{\prime}}{\sim
}X_{3}$, forming the cycle $(X_{1},Q_{5},X_{3},X_{2},X_{1})$. Hence $a\sim
x_{2}^{\prime}$. Applying Corollary \ref{Cor_C4}$(i)$ to $X_{3}\ $(and using
the private neighbour properties of $x_{1}^{\prime}$ and $x_{2}^{\prime}$), we
get $x_{1}^{\prime}\sim x_{2}^{\prime}$. Hence $x_{1}^{\prime},x_{2}^{\prime
},a$ have positive degree in $X_{3}$, and $X_{3}$-private neighbours
$x_{1},x_{2},x_{3}$, respectively, as asserted in (d). Thus, the flip-set of
$X_{3}$ using $x_{1},x_{2},x_{3}$ is $X$. Since $d_{T}(X,X_{3})=3$, Lemma
\ref{Lem_Not_Path} implies that $T$ contains $S(2,2)$.

$\bullet\hspace{0.1in}$Assume therefore that $a\notin\operatorname{EPN}%
(x_{3},X)$.

Similar to (b) we obtain that\ $a\ $is\ adjacent\ to$\ x_{1}$\ or$\ x_{2}%
,\ $but\ not\ to$\ x_{4},...,x_{r}$, and as in (c) we again obtain that $T$
contains $S(2,2)$. \medskip

\noindent\textbf{Case 2.2\hspace{0.1in} }Vertex $a$ is swapped into the
$\operatorname{IR}$-set in the second step. Then $X_{1}=\{x_{1}^{\prime}%
,x_{2},...,x_{r}\}$. Since $x_{2},x_{3}\nsim x_{1}^{\prime}$, Corollary
\ref{Cor_C4}$(i)$ applied to $X_{1}$ implies that $x_{2}\nsim x_{3}$. Since
$x_{1},x_{2},x_{3}$ have positive degree in $G[X]$, $x_{1}\sim x_{2},x_{3}$.

\noindent\underline{Suppose step 2 is $X_{1}\overset{x_{2}a}{\sim}X_{2}$.}

\noindent Then $X_{2}=\{x_{1}^{\prime},a,...,x_{r}\}$ and $X_{3}%
=\{x_{1}^{\prime},x_{2}^{\prime},...,x_{r}\}$. Now $X_{1}\overset{x_{2}%
x_{2}^{\prime}}{\sim}X_{3}$, a contradiction.

\noindent\underline{Suppose step 2 is $X_{1}\overset{x_{3}a}{\sim}X_{2}$.}

\noindent Then%
\[
X_{2}=\{x_{1}^{\prime},x_{2},a,x_{4},...,x_{r}\},\ X_{3}=\{x_{1}^{\prime
},x_{2}^{\prime},a,x_{4},...,x_{r}\}\text{\ and\ }X_{4}=\{x_{1}^{\prime}%
,x_{2}^{\prime},x_{3}^{\prime},x_{4},...,x_{r}\};
\]
also, $G[\{x_{3},x_{3}^{\prime},a\}]=K_{3}$. We show that

\begin{figure}[ptb]%
\centering
\includegraphics[
height=3.1721in,
width=5.1404in
]%
{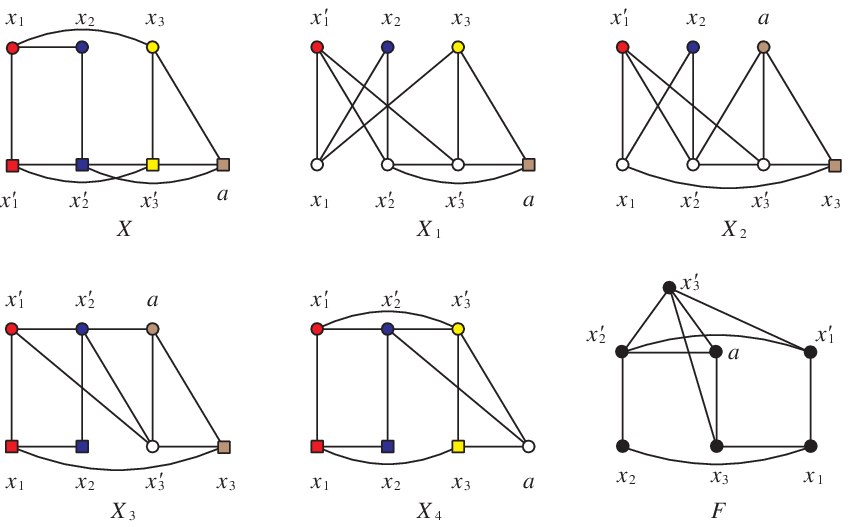}%
\caption{The $\operatorname{IR}$-sets $X$ and$\ X_{i},\ i=1,...,4$, in Case
2.2 when step 2 is $X_{1}\protect\overset{x_{3}a}{\sim}X_{2}$, with shared and
private neighbours, and the graph $F$}%
\label{Fig_552}%
\end{figure}

\begin{enumerate}
\item[(e)] $G[\{x_{1}^{\prime},x_{2}^{\prime},x_{3}^{\prime}\}]\cong K_{3}$,

\item[(f)] $X_{2}$ is independent, and

\item[(g)] $a\in\operatorname{EPN}(x_{3},X),$
\end{enumerate}

and deduce that $G[\{a,x_{1},x_{2},x_{3},x_{1}^{\prime},x_{2}^{\prime}%
,x_{3}^{\prime}\}]\cong F$, the graph in Figure \ref{Fig_Tree}, whose
$\operatorname{IR}$-graph is $\mathrm{Sp}(1,1;1,2)$. The sets $X,\ X_{i}%
,\ i=1,...,4$, and the graph $F$ are illustrated in Figure \ref{Fig_552}%
.\smallskip%
%TCIMACRO{\FRAME{ftbpFU}{5.1404in}{3.1721in}{0pt}{\Qcb{The $\operatorname{IR}%
%$-sets $X$ and$\ X_{i},\ i=1,...,4$, in Case 2.2 when step 2 is $X_{1}%
%\protect\overset{x_{3}a}{\sim}X_{2}$, with shared and private neighbours, and
%the graph $F$}}{\Qlb{Fig_552}}{thm_5_5_2.eps}%
%{\special{ language "Scientific Word";  type "GRAPHIC";
%maintain-aspect-ratio TRUE;  display "USEDEF";  valid_file "F";
%width 5.1404in;  height 3.1721in;  depth 0pt;  original-width 5.6818in;
%original-height 3.4938in;  cropleft "0";  croptop "1";  cropright "1";
%cropbottom "0";  filename 'Thm_5_5_2.eps';file-properties "XNPEU";}} }%
%BeginExpansion
%EndExpansion

\noindent(e)\hspace{0.1in}If $\{x_{1}^{\prime},x_{2}^{\prime},x_{3}\}$ is
independent, then $\{x_{1}^{\prime},x_{2}^{\prime},x_{3},...,x_{r}\}\neq
X_{3}$ is an $\operatorname{IR}(G)$-set adjacent to $X_{3}$ and $X_{4}$,
creating a cycle. Hence $x_{1}^{\prime}\sim x_{2}^{\prime}$.

If $x_{2}^{\prime}\nsim x_{3}^{\prime}$, then $\{x_{1},x_{2}^{\prime}%
,x_{3}^{\prime},x_{4},...,x_{r}\}\neq X_{3}$ is an $\operatorname{IR}(G)$-set
adjacent to $X_{4}$, also producing a contradiction. Hence $x_{2}^{\prime}\sim
x_{3}^{\prime}$. Similarly, if $x_{1}^{\prime}\nsim x_{3}^{\prime}$, then
$\{x_{1}^{\prime},x_{2},x_{3}^{\prime},x_{4},...,x_{r}\}\neq X_{3}$ is an
$\operatorname{IR}(G)$-set adjacent to $X_{4}$, hence $x_{1}^{\prime}\sim
x_{3}^{\prime}$. Therefore $G[\{x_{1}^{\prime},x_{2}^{\prime},x_{3}^{\prime
}\}]\cong K_{3}$ and (e) holds.\smallskip

\noindent(f)\hspace{0.1in}Applying Corollary \ref{Cor_C4}$(i)$ to $X_{2}$, we
see that either $x_{1}^{\prime}\sim a\sim x_{2}$ or $x_{1}^{\prime}\nsim
a\nsim x_{2}$. In the former case, let $u_{1},u_{2},u_{3}$ be the $X_{2}%
$-external private neighbours of $x_{1}^{\prime},x_{2},a$, and $Q_{6}$ the
resulting flip-set of $X_{2}$. Then $d_{T}(Q_{6},X_{2})\geq3$, a contradiction
by by Observation \ref{Ob_diam5}. Therefore $x_{1}^{\prime}\nsim a\nsim x_{2}%
$, that is, $\{x_{1}^{\prime},x_{2},a\}$ is independent.

Suppose $a$ is adjacent to a vertex in $\{x_{4},...,x_{r}\}$; say $a\sim
x_{4}$. Since $x_{1}^{\prime}\sim x_{2}^{\prime}$, we then see that
$x_{1}^{\prime},x_{2}^{\prime},a,x_{4}$ all have positive degree in $G[X_{3}%
]$, which is not the case (as $k=4$ was considered in Case 1). Therefore
$\{a,x_{4},...,x_{r}\}$ is independent, and since neither $x_{1}^{\prime}$ nor
$x_{2}$ is adjacent to any vertex in $\{a,x_{4},...,x_{r}\}$, $X_{2}$ is
independent, i.e., (f) holds. Noting that $x_{1}^{\prime}\sim x_{2}^{\prime}$
but $x_{1}^{\prime}\nsim a$, and applying Corollary \ref{Cor_C4}$(i)$ to
$X_{3}$, it follows that $x_{2}^{\prime}\sim a$.\smallskip\ 

\noindent(g)\hspace{0.1in}Suppose $a\notin\operatorname{EPN}(x_{3},X)$. Since
$X_{2}$ is independent, the only possibility is $a\sim x_{1}$. Note that
$x_{2}\in\operatorname{EPN}(x_{2}^{\prime},X_{3})$ and $x_{3}\in
\operatorname{EPN}(a,X_{3})$, but $\operatorname{EPN}(x_{1}^{\prime}%
,X_{3})\cap(X\cup\{x_{1}^{\prime},x_{2}^{\prime},x_{3}^{\prime}%
,a\})=\varnothing$. Since $x_{1}^{\prime}$ is not isolated in $G[X_{3}]$ (from
(e), it is adjacent to $x_{2}^{\prime}$), $\operatorname{EPN}(x_{1}^{\prime
},X_{3})\neq\varnothing$. Let $b\in\operatorname{EPN}(x_{1}^{\prime},X_{3})$,
$b\in V(G)-(\{a,x_{1}^{\prime},x_{2}^{\prime},x_{3}^{\prime}\}\cup X)$. Then
the flip-set $Q_{7}=\{b,x_{2},...,x_{r}\}$ of $X_{3}$ using $\{b,x_{2}%
,x_{3}\}$ is an $\operatorname{IR}(G)$-set. Now $Q_{7}\overset{bx_{1}}{\sim
}X_{1}$, and to avoid the triangle $(X_{0},X_{1},Q_{7},X_{0})$, $b\nsim x_{1}%
$. However, now $Q_{8}=\{x_{1},x_{2}^{\prime},b,x_{4},...,x_{r}\}$ is an
independent $\operatorname{IR}(G)$-set such that $d_{T}(X_{2},Q_{8})\geq3$, a
contradiction as before. This proves (g).\smallskip

Now $G[\{a,x_{1},x_{2},x_{3},x_{1}^{\prime},x_{2}^{\prime},x_{3}^{\prime
}\}]\cong F$, the graph in Figure \ref{Fig_Tree}, under the isomorphism
\[
a\rightarrow b,\ \ x_{1}\rightarrow g,\ \ x_{2}\rightarrow e,\ \ x_{3}%
\rightarrow f,\ \ x_{1}^{\prime}\rightarrow c,\ \ x_{2}^{\prime}\rightarrow
a,\ \ x_{3}^{\prime}\rightarrow d,
\]
as shown in Figure \ref{Fig_552}. Hence $T$ contains $\mathrm{Sp}(1,1;1,2)$.
Since $|V(T)|\leq7$, $T\cong\mathrm{Sp}(1,1;1,2)$.\medskip

\noindent\textbf{Case 2.3\hspace{0.1in}}Vertex $a$ is swapped third. Then
$X_{1}=\{x_{1}^{\prime},x_{2},...,x_{r}\}$, $X_{2}=\{x_{1}^{\prime}%
,x_{2}^{\prime},x_{3},...,x_{r}\}$ and $X_{3}=\{x_{1}^{\prime},x_{2}^{\prime
},a,...,x_{r}\}$. But then $X_{2}\sim_{T}X_{4}$, a contradiction.\medskip

This concludes the proofs of Case 2 and the theorem.~$\blacksquare$

\bigskip

The above results culminate in the following theorem on the realizability of
graphs as $\operatorname{IR}$-graphs.

\begin{theorem}
\label{Thm_Main}$(i)\hspace{0.1in}$Complete graphs are the only
$\operatorname{IR}$-graphs having universal vertices.

$(ii)\hspace{0.1in}$The cycles $C_{5},C_{6},C_{7}$ and the paths $P_{3}%
,P_{4},P_{5}$ are not $\operatorname{IR}$-graphs. 

$(iii)\hspace{0.1in}$The only connected $\operatorname{IR}$-graphs of order
$4$ are $K_{4}$ and $C_{4}$. 

$(iv)\hspace{0.1in}$The double star $S(2,2)$ is a smallest non-complete
$\operatorname{IR}$-tree.
\end{theorem}

\noindent\textbf{Proof.\hspace{0.1in}}Statement $(i)$ is Proposition
\ref{Prop_K1_n}. Considering cycles, an $\operatorname{IR}$-graph of diameter
$2$ has an induced $C_{4}$ (Proposition \ref{Prop_diam2}). Since
$\operatorname{diam}(C_{5})=2$, it follows that $C_{5}$ is not an
$\operatorname{IR}$-graph. Suppose there exists a graph $G$ such that
$G(\operatorname{IR})\cong C_{6}$ or $G(\operatorname{IR})\cong C_{7}$. By
Lemma \ref{Lem_Indep}, there exists an $\operatorname{IR}(G)$-set $X$ which is
not independent. By Corollary \ref{Cor_C4}, $G[X]$ contains at least two
edges, which means that $G[X]$ has at least three vertices of positive degree.
Since $\operatorname{diam}(C_{6})=\operatorname{diam}(C_{7})=3$, Corollary
\ref{Cor_Number_IR-sets1} states that $X$ has exactly three vertices of
positive degree. But then, by Lemma \ref{Lem_Not_Path}, $G(\operatorname{IR})$
contains an induced $4$-cycle or an induced double star $S(2,2)$, neither of
which is a subgraph of $C_{6}$ or $C_{7}$. We deduce that $C_{6}$ and $C_{7}$
are not $\operatorname{IR}$-graphs.

Since $\operatorname{diam}(P_{3})=2$, it is not an $\operatorname{IR}$-graph
by Proposition~\ref{Prop_diam2}. Similarly, by Proposition \ref{Prop_Diam3},
$P_{4}$ is not an $\operatorname{IR}$-graph, and by Theorem \ref{Thm_P5},
$P_{5}$ is not an $\operatorname{IR}$-graph. Hence $(ii)$ holds.

$(iii)\hspace{0.1in}$If $\operatorname{IR}(G)=1$, then $G$ is complete and
$G(\operatorname{IR})\cong G$; in particular, $K_{4}(\operatorname{IR})=K_{4}%
$. By $(i)$, the only non-complete connected $\operatorname{IR}$-graphs of
order $4$ have maximum degree $2$. By Proposition \ref{Prop_Qn}$(ii)$, $C_{4}$
is the $\operatorname{IR}$-graph of $2K_{2}$, and since $P_{4}$ is not an
$\operatorname{IR}$-graph, $(iii)$ is proved.

$(iv)\hspace{0.1in}$Since the only non-complete trees of order less than $5$
are $P_{3},P_{4}$ and $K_{1,3}$, which are not $\operatorname{IR}$-trees,
smallest non-complete $\operatorname{IR}$-trees have order at least $5$ and
diameter at least $3$. Since $S(2,2)$ is the unique smallest
$\operatorname{IR}$-tree with diameter~$3$, and $\mathrm{Sp}(1,1;1,2)$ is the
unique smallest $\operatorname{IR}$-tree with diameter~$4$ but has order $7$,
we consider trees of diameter at least $5$. All such trees have order at least
$7$, except for $P_{6}$, which has the same order as $S(2,2)$.~$\blacksquare$

\bigskip

We conjecture that $P_{6}$ is not an $\operatorname{IR}$-graph (see Conjecture
\ref{Con_paths}); if the conjecture is true, it will imply that $S(2,2)$ is
the unique smallest non-complete $\operatorname{IR}$-tree.

\section{Open Problems}

\label{Sec_Qs}A direct proof that $P_{5}$ is not an $\operatorname{IR}$-graph
is somewhat simpler than the proof of Theorem \ref{Thm_P5}, but not simple
enough to easily generalize to a proof that longer paths and cycles are not
$\operatorname{IR}$-graphs. Nevertheless, we believe this to be true and state
it as conjectures.

\begin{conjecture}
\label{Con_paths}$P_{n}$ is not an $\operatorname{IR}$-graph for each $n\geq3$.
\end{conjecture}

\begin{conjecture}
$C_{n}$ is not an $\operatorname{IR}$-graph for each $n\geq5$.
\end{conjecture}

\begin{problem}
Prove or disprove: Complete graphs and $K_{m}\boksie K_{n}$, where $m,n\geq2$,
are the only connected claw-free $\operatorname{IR}$-graphs.
\end{problem}

We showed in Lemma \ref{Lem_Indep} that if \emph{all} $\operatorname{IR}$-sets
of $G$ are independent, and the $\operatorname{IR}$-graph $H$ of $G$ is
connected and has order at least three, then $H$ contains a triangle or an
induced $C_{4}$. If $\operatorname{diam}(H)\geq3$, $H$ therefore contains a
vertex of degree at least three. The independent $\operatorname{IR}$-sets of
the graphs in Figures \ref{Fig_S22} and \ref{Fig_Tree}, which also have
non-independent $\operatorname{IR}$-sets, correspond to vertices of degree
three in their $\operatorname{IR}$-graphs. An affirmative answer to the next
question will be useful in proving that $C_{n}$ and $P_{n},n\geq6$, are not
$\operatorname{IR}$-trees.

\begin{problem}
Prove or disprove: if $G$ has an independent $\operatorname{IR}$-set and the
$\operatorname{IR}$-graph $H$ of $G$ is connected and has order at least
three, then $H$ has maximum degree at least three.
\end{problem}

\begin{conjecture}
\label{Con_doublestar}The only $\operatorname{IR}$-trees with diameter $3$ are
$S(2k,2k)$, where $k\geq1$.
\end{conjecture}

\begin{problem}
Determine which double spiders are $\operatorname{IR}$-trees.
\end{problem}

\begin{problem}
Characterise $\operatorname{IR}$-graphs having diameter $2$.
\end{problem}

\begin{problem}
As mentioned in Section \ref{Sec_Basic}, $K_{n}(\operatorname{IR})=K_{n}$.
Determine other graphs $G$ such that $G(\operatorname{IR})=G$, or show that
complete graphs are the only graphs with this property.\label{refs}
\end{problem}

\bigskip

\end{document}